\documentclass{article}
\usepackage{theorem}
\usepackage{amsmath}
\usepackage{amsfonts,amssymb} 
\usepackage{latexsym} 

\usepackage{xy}
\xyoption{all}

\IfFileExists{mathrsfs.sty}{\usepackage{mathrsfs}\let\mathcal\mathscr}{}

\newtheorem{theorem}{Theorem}[section]
\newtheorem{lemma}[theorem]{Lemma}
\newtheorem{prop}[theorem]{Proposition}
\newtheorem{corollary}[theorem]{Corollary}

{\theorembodyfont{\normalfont\rmfamily}
\newtheorem{definition}[theorem]{Definition}
\newtheorem{remark}[theorem]{Remark}
\newtheorem{question}[theorem]{Question}
}

\makeatletter
\def\ack{\vspace{.5\baselineskip}\noindent{\theorem@headerfont
Acknowledgement}\ \ }
    {\renewcommand{\theththm}{\ref{#1}}\begin{ththm}}
    {\end{ththm}}
\newtheorem{ththm}{Theorem}
\newenvironment{proof}[1][]%
{\def\proof@temp{#1}\par\noindent
\textsc{Proof}\ifx\proof@temp\@empty\else\
({\proof@temp})\fi\hspace{1em}}
{~\hfill{$\Box$}\par\vspace{.4\baselineskip}}
\def\operatorname#1{\mathop{\operator@font #1}\nolimits}%
\makeatother
\newcommand{\map}[1]{\stackrel{#1}{\longrightarrow}}
\newcommand{\CiM}{C^\infty(M)}
\newcommand{\bbnu}{[[\nu]]}

\newcommand{\A}{\mathbb{A}}
\newcommand{\TT}{\mathbb{T}}
\newcommand{\C}{\mathbb{C}}
\newcommand{\R}{\mathbb{R}}

\renewcommand{\H}{\mathcal{H}}

\newcommand{\g}{\mathfrak{g}}
\newcommand{\h}{\mathfrak{h}}
\newcommand{\Ad}{\operatorname{Ad}}
\newcommand{\ad}{\operatorname{ad}}

\newcommand{\Ker}{\operatorname{Ker}\,}
\newcommand{\Tr}{\operatorname{Tr}}

\newcommand{\Vol}{\operatorname{Vol}}
\newcommand{\End}{\operatorname{End}}

\newcommand{\id}{\operatorname{id}}
\newcommand{\Id}{\operatorname{Id}}

\newcommand{\half}{{\textstyle{\frac12}}}
\def\cyclic{\mathop{\kern0.9ex{{+}
\kern-2.2ex\raise-.28ex\hbox{\Large\hbox
{$\circlearrowright$}}}}\limits}

\newcommand{\suchthat}{\mathop{\,\vert\,}}
\newcommand{\jhrJ}{\mathcal{J}}
\newcommand{\jhrP}{\mathcal{P}}
\newcommand{\jhrV}{\mathcal{V}}
\newcommand{\jhrH}{\mathcal{H}}

\newcommand{\jhrm}{\mathfrak{m}}
\renewcommand{\sp}{\mathfrak{sp}}
\newcommand{\jhru}{\mathfrak{u}}

\newcommand{\widebar}[1]{{#1'}}
\newcommand{\nablabar}{\widebar{\nabla}}

\newcommand{\p}{{\mathfrak{p}}{}}

\newcommand{\s}{{\mathfrak{s}}{}}

\def\tref#1{Theorem~\ref{#1}}

\def\cref#1{Corollary~\ref{#1}}
\def\dref#1{Definition~\ref{#1}}

\newcommand{\Pf}{{\em Proof}. }
\newcommand{\EPf}
{%
\mbox{}%
\nolinebreak%
\hfill%
\rule{2mm}{2mm}%
\medbreak%
\par%
}

\newcommand{\trianglehorizontal}
{\mbox{\begin{picture}(40,40)\put(0,0){\line(1,0){20}}
\put(-1,-1.5){$\bullet$}
\put(37.5,-1.7){$\bullet$}\put(17.3,27.3){$\bullet$}
\put(40,0){\vector(-1,0){20}}\put(0,0){\vector(2,3){10}}
\put(10,15){\line(2,3){10}}\put(20,30){\vector(2,-3){10}}
\put(30,15){\line(2,-3){10}}\end{picture}}}
\newcommand{\Triangle}[3]{\mbox{\begin{picture}(60,100)
\put(10,20){$#1$}\put(60,20){$#3$}\put(35,60){$#2$}
\put(15,25){\begin{picture}(30,40)\put(0,0){\line(1,0){20}}
\put(-1,-1.5){$\bullet$}\put(37.5,-1.7){$\bullet$}
\put(17.3,27.3){$\bullet$}\put(40,0){\vector(-1,0){20}}
\put(0,0){\vector(2,3){10}}\put(10,15){\line(2,3){10}}
\put(20,30){\vector(2,-3){10}}\put(30,15){\line(2,-3){10}}
\end{picture}}\end{picture}}}

\newcommand{\Trianglefat}[3]{\mbox{\begin{picture}(60,100)
\put(10,20){$#1$}\put(60,20){$#3$}\put(35,60){$#2$}
\put(15,25){\begin{picture}(30,40)\thicklines

\put(0,0){\line(1,0){20}}
\put(-1,-1.5){$\bullet$}\put(37.5,-1.7){$\bullet$}
\put(17.3,27.3){$\bullet$}\put(40,0){\vector(-1,0){20}}
\put(0,0){\vector(2,3){10}}\put(10,15){\line(2,3){10}}
\put(20,30){\vector(2,-3){10}}\put(30,15){\line(2,-3){10}}
\end{picture}}\end{picture}}}

\newcommand{\carrehoriz}[1]{\mbox{\begin{picture}(100,100)
\put(10,10){\line(1,0){70}}\put(10,80){\line(1,0){70}}\put(10,10){\line(1,1){70}}

\put(10,80){\line(1,-1){70}}\put(30,30){\vector(1,1){2}}
\put(60,60){\vector(-1,-1){2}}\put(45,10){\vector(-1,0){2}}
\put(45,80){\vector(1,0){2}}\put(60,30){\vector(1,-1){2}}
\put(30,60){\vector(-1,1){2}}\put(8,8){$\bullet$}\put(8,78){$\bullet$}
\put(77,8){$\bullet$}\put(78,78){$\bullet$}\put(42,42){$\bullet$}\put(0,5){$d$}

\put(85,5){$c$}\put(85,80){$b$}\put(0,80){$a$}\put(50,43){#1}\end{picture}}}

\newcommand{\carrevertic}[1]{\mbox{\begin{picture}(100,100)
\put(10,10){\line(0,1){70}}\put(80,10){\line(0,1){70}}
\put(10,10){\line(1,1){70}}\put(10,80){\line(1,-1){70}}
\put(30,30){\vector(-1,-1){2}}\put(60,60){\vector(1,1){2}}
\put(10,45){\vector(0,1){2}}\put(80,45){\vector(0,-1){2}}
\put(60,30){\vector(-1,1){2}}\put(30,60){\vector(1,-1){2}}
\put(8,8){$\bullet$}\put(8,78){$\bullet$}\put(77,8){$\bullet$}
\put(78,78){$\bullet$}\put(42,42){$\bullet$}\put(0,5){$d$}\put(85,5){$c$}

\put(85,80){$b$}\put(0,80){$a$}\put(50,43){#1}\end{picture}}}
\newcommand{\carreverticfat}[1]{\mbox{\begin{picture}(100,100)
\thicklines
\put(10,10){\line(0,1){70}}\put(80,10){\line(0,1){70}}
\put(10,10){\line(1,1){70}}\put(10,80){\line(1,-1){70}}
\put(30,30){\vector(-1,-1){2}}\put(60,60){\vector(1,1){2}}
\put(10,45){\vector(0,1){2}}\put(80,45){\vector(0,-1){2}}
\put(60,30){\vector(-1,1){2}}\put(30,60){\vector(1,-1){2}}
\put(8,8){$\bullet$}\put(8,78){$\bullet$}\put(77,8){$\bullet$}
\put(78,78){$\bullet$}\put(42,42){$\bullet$}\put(0,5){$d$}\put(85,5){$c$}

\put(85,80){$b$}
\put(0,80){$a$}
\put(50,43){#1}
\end{picture}}}

\def\End{{\rm End}}
\def\mod{{\rm mod}} 
\def\ad{{\rm Ad}} 
\def\ad{{\rm ad}}

\def\ov{\overline}
\def\ot{\otimes}

\def\w{\wedge}
\def\ra{\rightarrow}

\newcommand{\SU}{{\mbox{\rm SU}}}

\newcommand{\G}{{\mbox{\rm G}}}
\newcommand{\T}{{\mbox{\rm T}}}

\newcommand{\eP}{{\mbox{\rm P}}}

\newcommand{\al}{\alpha}
\newcommand{\be}{\beta}

\newcommand{\om}{\omega}
\newcommand{\Om}{\Omega}
\renewcommand{\th}{\theta}

\renewcommand{\phi}{\varphi}
\def\rk{\mbox{\rm rk}}

\def\R{{\Bbb R}}
\def\C{{\Bbb C}}

\def\P{{\Bbb P}}

\def\su{{\frak {su}}}
\def\uu{{\frak {u}}}
\def\sl{{\frak {sl}}}
\def\sp{{\frak {sp}}}

\def\g{{\frak g}}
\def\h{{\frak h}}

\def\m{{\frak m}}

\def\s{{\frak s}}

\def\p{{\frak p}}

\renewcommand{\frak}{\mathfrak}
\renewcommand{\Bbb}{\mathbb}
\def\cC{{\cal C}}

\def\be{\begin{equation}}
\def\ee{\end{equation}}
\def\bi{\begin{enumerate}}
\def\ei{\end{enumerate}}
\def\ba{\begin{array}}
\def\ea{\end{array}}
\def\bea{\begin{eqnarray}}
\def\eea{\end{eqnarray}}
\def\ben{\begin{enumerate}}
\def\een{\end{enumerate}}

\def\ftnote#1{\def\footnotemark{}\footnote{#1}\setcounter{footnote}{0}}

\begin{document}

\title{Symplectic connections}

\author{
Pierre Bieliavsky${}^{1}$
\ftnote{\kern-4.5pt${}^{1}$Univ. Cath. Louvain,
D\'ept de Math, ch.\ du cyclotron 2, B-1348 Louvain-la-Neuve,
Belgium}
\\[20pt]
Michel Cahen${}^{2}$
\\[20pt]
Simone Gutt${}^{2,3}$
\ftnote{\kern-4.5pt${}^{2}$Universit\'e Libre
de Bruxelles, Campus Plaine, CP 218, B-1050~Brussels, Belgium}
\ftnote{\kern-4.5pt${}^{3}$Universit\'e de Metz, D\'ept. de Math., Ile du Saulcy,
F-57045~Metz Cedex 01, France}
\\[20pt]
John Rawnsley${}^{4}$ 
\ftnote{\kern-4.5pt${}^{4}$Mathematics Institute, University of Warwick,
Coventry\ \ CV4~7AL, United Kingdom}
\\[20pt]
Lorenz Schwachh\"ofer${}^{5}$
\ftnote{\kern-4.5pt${}^{5}$Math. Institut, Universit\"at 
Dortmund, Vogelpothsweg 87, D-44221~Dortmund, Germany}
\ftnote{\texttt{bieliavsky@math.ucl.ac.be, mcahen@ulb.ac.be, sgutt@ulb.ac.be,}}
\ftnote{\texttt{j.rawnsley@warwick.ac.uk, 
lschwach@math.uni-dortmund.de}}
 }
 
\date{}
\setcounter{page}{0}

\maketitle

\thispagestyle{empty}

\begin{abstract}
This article is an overview of the results obtained in recent years on
symplectic connections. We present what is known about preferred
connections (critical points of a variational principle). The class of
Ricci-type connections (for which the curvature is entirely determined
by the Ricci tensor) is described in detail, as well as its far reaching
generalization to special connections. A twistorial construction shows a
relation between Ricci-type connections and complex geometry. We give a
construction of Ricci-flat symplectic connections. We end up by presenting,
 through an explicit example, an approach
to noncommutative symplectic symmetric spaces.
\end{abstract}

\vspace{2cm}

\centerline{\small math.SG/0511194 v2, May 2006. Section 
\ref{subsect:twistorreduction} rewritten.}

\newpage
\pagenumbering{roman}
\tableofcontents
\enlargethispage{3\baselineskip}

\newpage
\pagenumbering{arabic}

\section{Introduction}

Symplectic geometry is by nature non-local. This is emphasized in 
particular by the classical Darboux theorem. The introduction
on a given symplectic manifold of a symplectic connection, which 
is a tool adapted to local computations, may seem inappropriate.
The aim of this survey is to show that symplectic geometry in the
presence of special symplectic connections becomes highly rigid.
More precisely there exists a family of universal models such
that any symplectic manifold admitting such a special symplectic 
connection is locally symplectically and affinely equivalent to a
particular model of the family. This local rigidity becomes a global 
rigidity if one requires compactness and simple connectedness.

It also appears that a twistor bundle over some of these
symplectic manifolds with special connections admits a
natural structure of complex analytic manifold.

Among the universal models there are certain symmetric symplectic 
spaces. These are particular manifolds where quantisation
(whether formal or convergent) may be performed explicitly;
hence they are a good framework for non-commutative geometry.

Although this survey is not exhaustive we have taken care to explain
various criteria for choosing particular symplectic connections.
Construction of quantisation based on a  choice of
symplectic connection is reduced here to one example, but we believe
this example shows possibilities for development.
We hope that this overview may lead others to investigate the
interplay of symplectic geometry and symplectic connections.

\section{Definitions  and basic facts about symplectic connections}

\subsection{Existence and non-uniqueness}
\begin{definition}
Let $(M,\omega)$ be a smooth symplectic  manifold of dimension $2n$
(i.e.  $\omega$ is a closed non-degenerate $2$-form on $M$). A
{\bf{symplectic connection}} on $(M,\omega)$ is a smooth linear
connection $\nabla$ such that:

-- its torsion $T^\nabla$ vanishes\\
 $~\qquad (\Leftrightarrow
T^\nabla(X,Y):= \nabla_XY-\nabla_YX-[X,Y]=0$);

-- the symplectic form $\omega$ is parallel \\
$~\qquad(\Leftrightarrow(\nabla_X\omega)(Y,Z):=X(\omega(Y,Z))
-\omega(\nabla_X Y,Z)-\omega(Y,\nabla_X Z)=0$).
\end{definition}

To prove the existence of such a connection,
take $\nabla^0$  any torsion free linear connection (for instance,
the Levi Civita connection associated to a metric $g$ on $M$).
Consider the tensor $N$ on $M$ defined by 
$$
\nabla^0_X\omega (Y,Z)=:\omega(N(X,Y),Z). 
$$
Since $\omega$ is skewsymmetric we have  
$\omega(N(X,Y),Z)=-\omega(N(X,Z),Y)$ and 
since $\omega$ is closed we have  $\cyclic_{XYZ}\omega(N(X,Y),Z)=0$ where
$\cyclic_{XYZ}$ denotes the sum over cyclic permutations of the indices
$X,Y$ and $Z$.
Define 
\[
\nabla_X Y:=\nabla^0_X Y + \frac{1}{3}N(X,Y)+\frac{1}{3}N(Y,X).
\] 
Then $\nabla$ is torsion free and:
\begin{eqnarray*}
 \nabla_X\omega (Y,Z)
&~&
=X(\omega(Y,Z))
-\omega(\nabla_X Y,Z)-\omega(Y,\nabla_X Z)
\\
&~&
=\nabla^0_X\omega (Y,Z)
-\frac{1}{3}\omega(N(X,Y),Z)-\frac{1}{3}\omega(N(Y,X),Z)
\\
&~&\quad
-\frac{1}{3}\omega(Y,N(X,Z))-\frac{1}{3}\omega(Y,N(Z,X))
\\
&~&
=\nabla^0_X\omega (Y,Z)-\frac{1}{3}\omega(N(X,Y),Z)
-\frac{1}{3}\omega(N(Y,X),Z)
\\
&~&\quad
-\frac{1}{3}\omega(N(X,Y),Z)
-\frac{1}{3}\omega(N(Z,Y),X)
\\
&~&
=(1-\frac{1}{3}-\frac{1}{3})\omega(N(X,Y),Z)+\frac{1}{3}\omega(N(X,Z),Y)=0
\end{eqnarray*}
so  the linear connection $\nabla$ is symplectic.

We shall now see how (non)-unique is a symplectic connection.
Take  $\nabla$ symplectic; then
$\nabla'_X Y:=\nabla_X Y + S(X,Y)$ is symplectic 
if and only if $S(X,Y)=S(Y,X)$ (torsion free) and
\begin{eqnarray*}
0&~&=\nabla'_X\omega(Y,Z)\\
&~&=\nabla_X\omega(Y,Z)-\omega(S(X,Y),Z)-\omega(Y,S(X,Z))\\
&~& =-\omega(S(X,Y),Z)+\omega(S(X,Z),Y),
\end{eqnarray*}
hence if and only if
$\omega(S(X,Y),Z)$ is totally symmetric.

Summarising we can now state the well known result:
\begin{theorem}
On a symplectic manifold $(M,\omega)$ there always exist symplectic connections. 
The set of symplectic connections
is an affine space modelled on the space of contravariant symmetric 
$3$-tensor fields on $M$, $\Gamma^\infty(S^3TM).$
\end{theorem}

\subsection{Where do symplectic connections arise?}

The notion of symplectic connection
is intimately related to that of natural formal
deformation quantisation at order $2$.
Quantisation of a classical system is a way to pass from classical
to quantum results.
Deformation quantisation was introduced by Flato,  Lichnerowicz
 and  Sternheimer in \cite{FLS} and  in \cite{BayenB}; they

``suggest that quantisation be understood as a deformation of the
structure of the algebra of classical observables rather than a radical change
in the nature of the observables.''

So deformation quantisation is  defined in terms
of a star product which is a formal
deformation of the algebraic structure of the space of smooth functions
on a symplectic (or more generally a Poisson) manifold. 
The associative structure given by the
usual product of functions and the Lie structure given by the Poisson
bracket are simultaneously deformed. Let us recall that 
if $(M,\omega)$ is a symplectic manifold and 
 if $u,v \in \CiM$,  the Poisson bracket of $u$ and $v$ is defined by
$$
\{u,v\}: = X_u(v) = \omega(X_v,X_u),
$$
where $X_u$ denotes the  Hamiltonian vector field 
corresponding to the function $u$,
i.e. such that $
i(X_u)\omega = du.$
\begin{definition} 
A \textbf{star product}
on a symplectic manifold $(M,\omega)$ is a bilinear map
\[
\CiM\times \CiM \to \CiM\bbnu ~~ (u,v) \mapsto  u*_\nu v: =
\sum_{r\ge0} \nu^rC_r(u,v)
\]
such that\\
$(u*v)*w = u*(v*w)$ (when extended $\R\bbnu$ linearly);\\
$C_0(u,v) = uv \qquad C_1(u,v)-C_1(v,u) = \{u,v\}$;\\
$1*u = u*1 = u$.\\ 
If all the $C_r$'s are bidifferential operators; one speaks of a 
{\bf differential star product}; if, furthermore, each
$C_r$ is of order $\le r$ in each argument, one speaks of a 
{\bf natural star product}.
\end{definition}
The link between symplectic connections and star products appear
already in the seminal paper \cite{BayenB} where the authors observe that
if there is a flat symplectic connection $\nabla$ on $(M,\omega)$,
one can generalise the classical formula for Moyal star product $*_M$ 
defined on $\R^{2n}$
with a constant symplectic $2$-form. 

Fedosov, proved more generally that given any symplectic connection
$\nabla$, one can construct a star product (in \cite{GRS} it was
proposed that a triple $(M,\omega,\nabla)$ be known as a Fedosov
manifold):

\begin{theorem} \cite{Fed}
Given a symplectic connection $\nabla$ 
 and a sequence $\Omega=\sum_{k=1}^\infty \nu^k\omega_k$ of closed $2$-forms 
 on a symplectic manifold $(M,\omega)$, 
 one can build a star product $*_{\nabla,\Omega}$ on it. This is 
obtained by identifying the space $\CiM[[\nu]]$ with a subalgebra of 
the algebra of sections of a bundle of associative algebras
(called the Weyl bundle) on $M$. The subalgebra is the one of flat sections
of the  Weyl bundle, when this bundle is endowed with a flat connection whose
construction is determined by the  choices made of the connection on $M$ and
 of the sequence of closed 
$2$-forms on $M$. 
\end{theorem}

Reciprocally a natural star product determines a symplectic connection.
This was first observed by Lichnerowicz \cite{Lichne} for a restricted class of
star products.
\begin{theorem}{\cite{GRnatural}} A natural star product 
at order 2 determines a unique symplectic connection. 
\end{theorem}

\subsection{When is there a ``natural'' unique symplectic connection?}

To have a canonical choice of  symplectic 
connection on $(M,\omega)$,
one needs some extra structure on the manifold.

{\bf $\bullet$ Example 1: pseudo-K\"ahler manifolds}

\noindent Choose an almost complex  structure $J$ on $(M,\omega)$
[i.e. $J:TM\rightarrow TM$ is a bundle endomorphism so that
$J^2=-\Id$] so that $\omega(JX,JY)=\omega(X,Y)]$.

\noindent A symplectic connection $\nabla$ 
 preserves $J$ [i.e.  $\nabla J=0$]
if and only if it is the Levi Civita connection associated to the
pseudo Riemannian metric $g(X,Y):=\omega(X,JY)$. It is 
thus unique
and it only exists in a  
(pseudo-)K\"ahler situation. 

{\bf $\bullet$ Example 2: symmetric symplectic spaces}

\noindent Intuitively, a symmetric symplectic space is a  symplectic manifold
with symmetries attached to each of its points. Precisely:
\begin{definition}\label{def:symsymspaces}
A {\bf{symmetric symplectic space}} is a triple $(M,\omega,S)$ 
where $(M,\omega)$ is a symplectic manifold
and where $S$ is a smooth map $S:M\times M\rightarrow M$
such that, defining for any point $x\in M$ the map (called the symmetry at $x$):
$$
s_x:=S(x,\cdot):M\rightarrow M,
$$
each $s_x$ squares to the identity [$s_x^2=\Id$] and
is a symplectomorphism of $(M,\omega)$ [$s_x^*\omega=\omega$],\\
$x$ is an isolated fixed point of $s_x$, and
$s_xs_ys_x=s_{s_xy}$ for any $x,y\in M$.
\end{definition}
\begin{prop}\cite{Biel}
On a symmetric symplectic space, there is a unique  symplectic connection 
for which each $s_x$ is an affinity.
It is explicitly given by
$$
\omega_x(\nabla_XY,Z)=\frac{1}{2}X_x\omega(Y+{s_{x}}_\star Y,Z).
$$
\end{prop}

{\bf $\bullet$} In the two examples above, the choice of symplectic
connection was imposed by the presence of an additional structure. 
To select a ``small'' class of symplectic connections on a symplectic
manifold without any additional structure, one has to choose some
restrictive conditions. One way to proceed is to impose some system 
of equations on the curvature tensor.

\subsection{Curvature tensor of a  symplectic connection} 

The {\bf curvature tensor} $R^\nabla$ of a linear connection $\nabla$
is the $2$-form on $M$ with values in the endomorphisms of
the tangent bundle defined by
\begin{equation}
R^\nabla(X,Y)Z = \left(\nabla_X\nabla_Y - \nabla_Y\nabla_X -
\nabla_{[X,Y]}\right)Z
\end{equation}
for vector fields $X,Y,Z$ on $M$. If $\nabla$ is symplectic,
$R^\nabla_x(X,Y)$ has values in the symplectic Lie algebra
$sp(T_xM,\omega_x)=\{ A \in \End (T_xM)~\vert ~
\omega_x(Au,v)+\omega_x(u,Av)=0,~ \forall u,v \in T_xM\}$.

The  curvature tensor
 satisfies the first Bianchi identity
\[
\cyclic_{X,Y,Z} R^\nabla(X,Y)Z = 0
\]
where $\cyclic$ denotes the sum over the cyclic permutations of the listed
set of elements,and the second Bianchi identity
\[
\cyclic_{X,Y,Z} \left(\nabla_XR^\nabla\right)(Y,Z) = 0.
\]

The {\bf Ricci tensor} $r^\nabla$ is the  $2$-tensor
\begin{equation}
r^\nabla(X,Y) = \Tr\left( Z \mapsto R^\nabla(X,Z)Y\right).
\end{equation}
The first Bianchi identity implies that $r^\nabla$ is symmetric.

One can define a second trace
$r'_x(X,Y):=\sum_i\omega (R^\nabla_x(e_i,e^i)X,Y)$
where the $e_i$ constitute a basis of $T_xM$
and the $e^i$ constitute the dual basis of $T_xM$
(i.e. such that $\omega (e_i,e^j)=\delta_i^j)$.
Then Bianchi's first identity
implies that $r'=-2r^\nabla$.

Since the Ricci tensor is symmetric and one only has 
a skewsymmetric contravariant $2$-tensor on $M$
(the Poisson tensor related to the symplectic form)
there is {\bf no ``scalar curvature''}.

The {\bf symplectic curvature tensor} is defined as
\begin{equation}
{\underline{R}}^\nabla(X,Y,Z,T):= \omega(R^\nabla(X,Y)Z,T).
\end{equation}
It is antisymmetric in its first two arguments and symmetric in its last
two. Hence ${\underline{R}}_x^\nabla$ is in $\Lambda^2(T^*_xM)\otimes
S^2(T^*_xM)$. To understand  Bianchi's first identity, we introduce the
operators of symmetrisation and skewsymmetrisation arising in the Koszul
long exact sequence.

\subsubsection{The Koszul long exact sequence}

Given any finite dimensional vector space $V$, the Koszul long exact sequence
has the following form:
$$
0 \longrightarrow S^q(V) {\stackrel{a}{\longrightarrow}} V\otimes S^{q-1}(V)
{\stackrel{a}{\longrightarrow}} \Lambda^2 V\otimes S^{q-2}(V)
{\stackrel{a}{\longrightarrow}}\cdots{\stackrel{a}{\longrightarrow}}
\Lambda^{q-1}(V)\otimes V {\stackrel{a}{\longrightarrow}}
\Lambda^q(V)\longrightarrow 0
$$
where $a$ is the skewsymmetrisation operator:
$$
a(v^1\wedge\ldots\wedge v^q\otimes w^1\cdots w^p)=
\sum_{i=1}^pv^1\wedge\ldots\wedge v^q\wedge w^i\otimes w^1\cdots
w^{i-1}w^{i+1}\cdots w^p.
$$
The symmetrisation operator reads:
$$
s(v^1\wedge\ldots\wedge v^q\otimes w^1\cdots w^p)
\sum_{i=1}^q(-1)^{q-i}v^1\wedge\ldots v^{i-1}\wedge v^{i+1}\ldots\wedge v^q 
\otimes v^i\cdot w^1\cdots w^p.
$$
These two operators satisfy $a^2=0,~s^2=0,~(a_\circ s+s_\circ a)_{\vert_{\Lambda^q
V\otimes S^{p}(V)}}=(p+q)\Id.$

The first  Bianchi identity on the value at a point $x$ of the
symplectic curvature tensor takes the form:
$$
\cyclic_{X,Y,Z} {\underline{R}}^\nabla_x(X,Y,Z,T)= 0
\Leftrightarrow {\underline{R}}^\nabla_x \in \ker a \subset
 \Lambda^2(T^*_xM)\otimes S^2(T^*_xM).
 $$

The space ${\underline{\mathcal{R}}}_x$ of 4-tensors satisfying the
algebraic identities of a symplectic curvature tensor at $x$ is:
 $$
{\underline{\mathcal{R}}}_x:=\ker a_{\vert_{
 \Lambda^2(V)\otimes S^2(V)}}\simeq \left( V\otimes S^3(V)\right)/S^4(V)
 \qquad {\mathrm{for~}}~ V=T^*_xM.
 $$

\subsubsection{Decomposition of the curvature}

The group $Sp(T_xM, \omega_x)=
\{~A\in\End(T_xM)~\vert~\omega_x(Au,Av)=\omega_x(u,v)~\forall u,v\in T_xM~\}$ 
acts on $V=T^*_xM$ and thus  on  ${\underline{\mathcal{R}}}_x\simeq
\left( V\otimes S^3(V)\right)/S^4(V)$. Under this action the space $
V\otimes S^3(V)$, in dimension $2n \geq 4$, decomposes into three
irreducible subspaces ($ S^4(V)\oplus S'^2(V)\oplus W$ where $S'^2(V)=
a(s(\omega_x\otimes S^2(V)))\sim S^2(V) $so that:
$$
{\underline{\mathcal{R}}}_x = {\underline{{\mathcal E}}}_x 
\oplus {\underline{{\mathcal W}}}_x.
$$
The decomposition of the symplectic curvature tensor
${\underline{R}}^\nabla_x$ into its ${\underline{{\mathcal E}}}_x$
component (denoted ${\underline{E}}^\nabla_x$) and its
${\underline{{\mathcal W}}}_x$ component (denoted
${\underline{W}}^\nabla_x$) ,
$$
{\underline{R}}^\nabla_x = {\underline{E}}^\nabla_x +{\underline{W}}^\nabla_x,
$$ is given by 
\begin{eqnarray*}
{\underline{E}}^\nabla_x(X,Y,Z,T) &=& -\frac{1}{2(n+1)} \big[ 2 \omega_x
(X,Y) r^\nabla_x (Z,T) + \omega_x(X,Z) r^\nabla_x (Y,T)\\[2mm] &\,& +
\omega_x (X,T) r^\nabla_x (Y,Z) - \omega_x(Y,Z) r^\nabla_x(X,T) -
\omega_x (Y,T) r^\nabla_x(X,Z)\big].
\end{eqnarray*}
The corresponding decomposition of the curvature tensor (see Vaisman  
\cite{Vaisman}) has the form
\begin{equation}
{{R}}^\nabla_x = {{E}}^\nabla_x +{{W}}^\nabla_x,
\end{equation}
where
\begin{eqnarray}
E^\nabla(X,Y)Z &=& 
{\textstyle{\frac1{2n+2}}}\biggl(2\omega(X,Y)\rho^\nabla Z
 + \omega(X,Z)\rho^\nabla Y - \omega(Y,Z)\rho^\nabla X\\ \nonumber
&~&\quad\quad+ \omega(X,\rho^\nabla Z)Y - \omega(Y,\rho^\nabla Z)X
\biggr)
\end{eqnarray}
with $r^\nabla$ converted into an endomorphism $\rho^\nabla$ given by
\begin{equation}\label{eq:ricciend}
\omega(X, \rho^\nabla Y) = r^\nabla(X,Y).
\end{equation}
\begin{definition}\label{def:Riccitype}
A symplectic connection $\nabla$ on $(M,\omega)$ will be said to be
\textbf{of Ricci-type} if $W^\nabla=0$; it will be said to be 
\textbf{Ricci-flat} if $E^\nabla=0$ (hence if and only if $r^\nabla=0$).
\end{definition}
One can combine  restricting the holonomy algebra $\g\subset
sp(\R^{2n},\Omega)$ and the vanishing of some components of the
curvature when the curvature is decomposed into
irreducible components under the action of $\g$ to define 
special symplectic connections; these will appear  in section
\ref{section:SpecialSymplectic}.

\subsection{Variational principle}

To select symplectic connections through a variational principle
\cite{BourgeoisCahen}, one can consider a Lagrangian $L(R^\nabla)$,
which is a polynomial in the curvature of the connection $\nabla$,
invariant under the action of the symplectic group
$$
\int_M L(R^\nabla)\omega^n.
$$
\, From what we have seen before, there is no invariant polynomial
of degree $1$ in the curvature, so the easiest choice is
a polynomial  of  degree $2$ in $R^\nabla$. The space of
degree $2$ polynomials in the curvature which are invariant under
the action of the symplectic group
is $2$-dimensional and spanned by
$E^\nabla\cdot E^\nabla$
and $W^\nabla\cdot W^\nabla$ (or, equivalently by $R^\nabla\cdot
R^\nabla$ and $r^\nabla\cdot r^\nabla$) where $\cdot$ denotes the
symmetric function-valued product of tensors induced by $\omega$ and
$S \cdot T$, for $S$ and $T$  tensor-fields on $M$ of the same type, 
is given in local coordinates by
$$
S \cdot T=(\omega^{-1})^{i_1i'_1}\cdots(\omega^{-1})^{i_pi'_p}
\omega_{j_1j'_1}\cdots \omega_{j_qj'_q}
S_{i_1\ldots i_p}^{j_1\ldots j_q}T_{i'_1\ldots i'_p}^{j'_1\ldots j'_q}.
$$

From Chern--Weil theory we know that the first Pontryagin class is
represented by a $4$-form $P_1(\nabla)$ which is built from an invariant
combination of the curvature and the symplectic form. In fact
$$
P_1(\nabla)\wedge \omega^{n-2}= \frac{1}{16 \pi^2}[r^\nabla\cdot
r^\nabla-\half R^\nabla\cdot R^\nabla] \omega^n
$$
Since this combination will be constant under variations, all
non-trivial Euler equations coming from a variational principle built
from a second order invariant polynomial in the curvature are the same
and take the form:
\begin{equation}\label{eq:preferred}
 \cyclic_{X,Y,Z} (\nabla_X r^\nabla)(Y,Z)=0.
\end{equation}
 \begin{definition} A symplectic connection $\nabla$ is said to be
 {\bf{preferred}} if it is a solution of Equation
 (\ref{eq:preferred}).
\end{definition}

\section{Preferred symplectic connections}

The preferred symplectic connections $\nabla$ on a symplectic manifold 
$(M,\omega)$ of dimension $2n$ are critical points of the functional
$$
\int_M \Tr(\rho^\nabla)^2\frac{\omega^n}{n!}
$$
where $\rho^\nabla$ is the Ricci endomorphism 
as previously defined in (\ref{eq:ricciend}).
They obey the system
of second order partial differential equations
$$
\cyclic_{X,Y,Z} (\nabla_X r^\nabla)(Y,Z)=0
$$
where $r^\nabla$ is the Ricci tensor of the connection $\nabla$.\\ 
The basic problem is to determine if on a given manifold $(M,\omega)$
there exists a preferred connection and ``how many'' different ones may
occur; two solutions are different in this context if they are not
related through a symplectic diffeomorphism.

This basic problem is essentially solved when $(M,\omega)$
is a compact orientable surface. We give also a partial answer for
the standard symplectic plane $(\R^2,\Omega_0)$.

 When the dimension of $M$ is $\ge 4$, we can give large classes of
 examples of preferred connections:\\
(i) examples of symplectic manifolds admitting a Ricci-flat connection;\\
(ii) examples of symplectic manifolds admitting a ``Ricci-parallel''
connection;\\
(iii) examples of homogeneous symplectic manifolds admitting a
homogeneous preferred symplectic connection.

\subsection{Preferred symplectic connections in dimension $2$}

The study of preferred connections in dimension $2$ relies on a function
$\beta$, depending on the connection, which first appears in
\begin{lemma}\cite{BourgeoisCahen}
Let $(M,\omega)$ be a symplectic surface and let $\nabla$
be a preferred symplectic connection on it. Then 
\begin{itemize} 
\item[(i)] there exists a $1$-form $u$ such that, for any vector fields
$X,Y,Z$ one has
\begin{equation}\label{u2}
(\nabla_X r^\nabla)(Y,Z)= \omega(Y,X)u(Z)+\omega(Z,X)u(Y);
\end{equation}
\item[(ii)] there exists a function $\beta$ such that
\begin{equation}\label{beta}
\nabla_X u = \beta \omega;
\end{equation}
\item[(iii)] defining the vector field ${\bar{u}}$ by
$$
i({\bar{u}})\omega=u,
$$
one has, for any vector fields $X,Y$:
\begin{eqnarray*}
X\beta&=&-r^\nabla(X,{\bar{u}}),\cr
(\nabla^2\beta)(X,Y)&=&-u(X)u(Y)+\beta r^\nabla(X,Y);
\end{eqnarray*}
\item[(iv)] there exist two real numbers $A$ and $B$ such that
\begin{eqnarray*}
r^\nabla({\bar{u}},{\bar{u}})&=&\beta^2+B,\cr
\frac{1}{4}\Tr (\rho^\nabla)^2&=&\beta+A.
\end{eqnarray*}
\end{itemize}
\end{lemma}

Hence $\nabla$ is locally symmetric (i.e. $\nabla R^\nabla=0$) if and
only if $\beta=0$. If $M$ is compact, the second property above shows
that $\beta$ can not be a non-vanishing constant (it would indeed imply
that the symplectic $2$-form $\omega$ be exact). If $\beta$ is not a
constant and $M$ is compact, a detailed study of the critical points of
$\beta$ permits to show that no compact symplectic surface of positive
genus admits a non-locally symmetric preferred symplectic connection.
The case of the $2$-sphere is more delicate and requires precise
estimates. The conclusion is

\begin{theorem}\cite{BourgeoisCahen} 
A preferred connection on a compact symplectic surface is necessarily
locally symmetric.
\end{theorem}

Globally symmetric symplectic surfaces may be completely described. 
Locally symmetric symplectic surfaces may be simply related to the
globally symmetric ones provided the connection is geodesically
complete. A case by case analysis leads to

\begin{theorem}\cite{BourgeoisCahen,Biel}
Let $(M,\omega,\nabla)$ be a compact symplectic surface endowed with a
complete locally symmetric symplectic connection. Then, up to
diffeomorphisms, either
\begin{itemize}
\item $M$ is the $2$-sphere $S^2$ with $\omega$ a multiple of the
standard volume form  and with $\nabla$ the Levi Civita connection of
the standard Riemannian metric with constant positive  curvature equal
to $1$;
\item $M$ is the torus $T^2$ with $\omega$ a multiple of the standard
invariant volume form and with $\nabla$ a  flat affine symplectic
connection;
\item $M$ is a surface $\Sigma_g$ of genus $g\ge 2$ with $\nabla$
a multiple of the standard volume form inherited from the disk
and with $\nabla$ the connection associated to a metric
$h$ of constant negative curvature equal to $-1$.
\end{itemize}
\end{theorem}

What can be said about the existence and the number of different
complete preferred connections on a given compact symplectic surface?

Consider the case of $S^2$ with a symplectic structure $\omega$. There
exists a positive real number $k$ so that $\int_{S^2}\omega=k
\int_{S^2}\omega_0$ where $\omega_0$ is the standard symplectic
structure on $S^2$ defining the same orientation as $\omega$. Thus the
de Rham cohomology classes defined by $\omega$ and $k\omega_0$ are the
same. Furthermore $(1-t)\omega+ tk\omega_0$ is symplectic for any $t\in
[\,0\,,\,1\,]$ and defines the same cohomology class. By Moser's
stability theorem (see \cite{McDuffSal}), there exists a $1$-parametric
family $\phi_t$ of diffeomorphisms of the sphere $S^2$ such that
$\phi_0=\id$ and $\phi_1^*\omega=k\omega_0$. If $g_0$ is the standard
metric on $S^2$, define the metric $g$ on the sphere to be such that
$\phi_1^*g=kg_0$. The Levi Civita connection associated to $g$ is 
clearly symplectic relative to $\omega$ and symmetric, hence preferred.
Furthermore, if $\nabla'$ is another symplectic connection on
$(S^2,\omega)$ which is preferred and complete, it is automatically
symmetric since $S^2$ is simply connected and it coincides with
$\nabla$. Thus

\begin{theorem}\cite{BourgeoisCahen}
On $(S^2,\omega)$ there exists a complete symplectic preferred
connection and any complete preferred connection is the image of that
one through a symplectomorphism.
\end{theorem}
The result on the existence is obtained similarly on the torus $T^2$ or
the surface $\Sigma_g$ endowed with any symplectic structure. In these
cases, there is no unicity (because of the choice of the affine
connection or the choice of the metric).

The non-compact situation is more complicated. One can show that on the
plane endowed with the standard constant symplectic structure
$(\R^2,\Omega_0)$, there exist five affinely distinct globally symmetric
complete symplectic connections. There also exist a $2$-parametric
family of non-homogeneous preferred symplectic connections which are not
locally symmetric.

\subsection{Ricci-flat connections}
 
A symplectic connection is said to be Ricci-flat if its Ricci tensor
$r^\nabla$ vanishes. Obviously, those give examples of preferred
connections! We explain in Section \ref{sec:Ricciflat} a construction of
Ricci-flat connections: when a symplectic manifold $(M,\omega)$ of
dimension $2n\ge 4$ is the first element of a contact quadruple, any
symplectic connection $\nabla$ on $(M,\omega)$ can be lifted to define a
Ricci-flat symplectic connection on a certain symplectic manifold
$(P,\omega')$ of dimension $2n+2$. This procedure gives examples of
Ricci-flat, non-flat, symplectic connections in any dimension $\ge 6$.

\subsection{Ricci-parallel symplectic  connections}
 
Let $(M,\omega,\nabla)$ be a symplectic manifold with a symplectic
connection; assume there exists on this manifold a  compatible almost
complex structure $J$ (i.e. at each point $x$ in $M$, $J_x\in
\End(T_xM)\,,~J_x^2=-\Id\,,~ \omega(JX,JY)=\omega(X,Y)$) which is
parallel (i.e. $\nabla J=0$). Then the pseudo-riemannian metric $g$
defined by $g(X,Y):=\omega(X,JY)$ is also parallel, so the connection is
the Levi Civita connection for $g$ and the manifold  is pseudo-K\"ahler.
If this connection is preferred one has:

\begin{theorem}\cite{bib:CGR1} 
Let $(M,\omega,J)$ be a pseudo-K\"ahler manifold. If the
Levi Civita connection is preferred, then the Ricci tensor is parallel.
\end{theorem}

\noindent This gives examples of Ricci-parallel symplectic connections.

Consider now a symplectic manifold $(M,\omega)$ endowed with a
Ricci-parallel symplectic connection $\nabla$. Write
$$
T_xM^{\C}=\oplus_{\lambda\in spec}(T_xM^{\C})_\lambda
$$
where $(T_xM^{\C})_\lambda$ is the generalised eigenspace for
$\rho_x^\nabla$,i.e.
$$
(T_xM^{\C})_\lambda=\{X\in T_xM^{\C}\,\vert\,(\rho_x^\nabla-\lambda)^{2n}X=0\}
$$
and $spec$ denotes the set of different eigenvalues of $\rho_x^\nabla$ on
$T_xM^{\C}$.

Observe that  $(T_xM^{\C})_\lambda$ and $(T_xM^{\C})_\mu$ are orthogonal
with respect to the symplectic form $\omega_x$ unless $\mu=-\lambda$.
We define real symplectic subspaces of $T_xM$:\\
$R_\lambda$, for each real positive eigenvalue $\lambda$, so that the 
complexification of $R_\lambda$ is given by $R_\lambda^\C=
(T_xM^{\C})_\lambda\oplus (T_xM^{\C})_{-\lambda}$; we denote by $R$ the
set of such eigenvalues;\\
$I_\lambda$ for each  purely imaginary eigenvalue $\lambda=ia,\, a>0$,
so that its complexification is $I_\lambda^\C= (T_xM^{\C})_\lambda\oplus
(T_xM^{\C})_{-\lambda}$; we denote by $I$
the set of such eigenvalues;\\
$C_\lambda$ for each  eigenvalue $\lambda=a+ib,\, a>0,b>0$, so that its
complexification is $C_\lambda^\C= (T_xM^{\C})_\lambda\oplus
(T_xM^{\C})_{-\lambda}\oplus (T_xM^{\C})_{\overline{\lambda}}\oplus
(T_xM^{\C})_{-{\overline{\lambda}}}$; we denote by $C$ the set of such
eigenvalues.\\
They give a symplectic orthogonal decomposition of $T_xM$:
$$
T_xM=T_xM_0\oplus(\oplus_{\lambda\in R}R_\lambda)\oplus
(\oplus_{\lambda\in I}I_\lambda)\oplus(\oplus_{\lambda\in C}C_\lambda).
$$

Ricci-parallel implies that the Ricci endomorphism $\rho^\nabla$ 
commutes with all curvature endomorphisms so that all 
distributions corresponding to the above defined  subspaces
are parallel. We get

\begin{theorem}\cite{bib:CGR1,Boubel}
Let $(M,\omega)$ be a $2n$-dimensional symplectic manifold endowed with
a symplectic connection $\nabla$ whose Ricci tensor is parallel. Assume
that the Ricci tensor is not degenerate. Then
\begin{itemize}
\item[1] the connection $\nabla$ is the Levi Civita connection
associated to the metric defined by the Ricci tensor $r^\nabla$;
\item[2] the distributions $R_\lambda,I_\lambda,C_\lambda$ are parallel,
symplectic, and $\rho^\nabla$ restricted to any of these is semisimple;
\item[3] if $M$ is simply connected and $\nabla$ complete, then $M$ is
symplectomorphic and affinely equivalent to the product of the
symplectic submanifolds corresponding to the integral leaves of the
distributions $R_\lambda,I_\lambda,C_\lambda$;
\item[4] the manifolds corresponding to the leaves of $R_\lambda$
admit two parallel, transverse Lagrangian foliations; those corresponding
to $I_\lambda$ are K\"ahler-Einstein manifolds; those corresponding to 
$C_\lambda$ have dimension $4k$ if $k$ is the multiplicity of
$\lambda$ and also admit two parallel transverse Lagrangian foliations;
\item[5] if all the eigenvalues of $\rho^\nabla$ have multiplicity one,
the factors are two dimensional or four dimensional symplectic
symmetric spaces.
\end{itemize}
\end{theorem}

\subsection{Homogeneous preferred connections}

If $(M,\omega)$ is a connected, simply connected, compact, homogeneous,
symplectic manifold, a classical result of Kostant tells us that 
$(M,\omega)$ is symplectomorphic to a coadjoint orbit of a compact
semi simple Lie group $G$, endowed with its standard 
Lie-Kirillov-Kostant-Souriau symplectic structure.

Let $\nabla$ be a symplectic connection on $(M,\omega)$ stable by
the action of $G$. Let $p\in M$ and let $\pi:G\rightarrow M~g\mapsto g.p$
be the canonical projection related to the choice of the base point $p$.
Let $H$ be the stabilizer of $p$ in $G$ 
($H:=\{\,g\in GÖ,\vert\, g.p=p\,\}$) and let $\g$ (resp. $\h$) be
the Lie algebra of $G$ (resp. $H$). For any element $X\in\g$, let $X^*$
denotes the fundamental vector field on $M$ associated to $X$, 
i.e. $X^*_x=\frac{d}{dt\vert_0} \exp -tX.x$.
Write $\g=\h\oplus {\mathfrak{m}}$ where ${\mathfrak{m}}$ is the orthogonal
to $\h$ relative to the Killing form.
Each tangent vector at $p$ is the value at $p$ of a
(unique) fundamental vector field $Y^*$ with $Y\in {\mathfrak{m}}$.
Define a map $D:{\mathfrak{m}}\rightarrow \End({\mathfrak{m}})$ by
$$
(\nabla_{X^*}{Y^*})_{p}=(D(X)Y)^*_p.
$$
This map $D$ entirely determines the $G$-invariant connection $\nabla$.
It satisfies two conditions:
\begin{eqnarray}\label{D}
D(X)Y-D(Y)X-\pi_{\mathfrak{m}}([X,Y])&=&0\cr
\Omega(D(Y)X,Z)+\Omega(Y,D(Z)X)&=&0
\end{eqnarray}
where $\pi_{\mathfrak{m}}$ is the projection of $\g$ on ${\mathfrak{m}}$
 relatively to the decomposition
$\g=\h\oplus {\mathfrak{m}}$, and where 
$\Omega=(\pi^*\omega_p)_{\vert_{{\mathfrak{m}}\times{\mathfrak{m}}}}$.
The space of invariant homogeneous symplectic connections on $M$
may be identified to the space of maps 
$D:{\mathfrak{m}}\rightarrow \End({\mathfrak{m}})$ 
satisfying the conditions \ref{D}. Such a connection is preferred
if it is a critical point of the functional
$$
{\mathcal{I}}(\nabla)=\int_M \Tr(\rho^\nabla)^2\frac{\omega^n}{n!}.
$$
By Palais's principle, to determine the $G$-invariant critical points of 
${\mathcal{I}}$, it is sufficient to determine the critical points
of the restriction of ${\mathcal{I}}$ to the space of $G$-invariant
connections \cite{Palais}. For an invariant connection, ${\mathcal{I}}$
reduces to $\Tr(\rho^\nabla)^2 \Vol (M)$ and one shows that the Ricci 
tensor is a polynomial of degree $2$ on $D$.
\begin{lemma}\label{lemmapol1}\cite{CGR2} 
The functional ${\mathcal{I}}$ is a fourth order polynomial on $D$.
\end{lemma}
Using the structure of semisimple algebras of compact type, one shows 
that this polynomial is non-negative and that the homogeneous polynomial 
corresponding to the terms of degree $4$ is strictly positive outside 
the origin. In fact one replaces $D$ by 
$D'=D-\half \pi_{{\mathfrak{m}}}\circ \ad$; the conditions \ref{D} simply
express the fact that the $3$-form $\Omega(D'(\cdot)\cdot,\cdot)$ is
completely symmetric.
\begin{lemma}\label{lemmapol2}\cite{CGR2} 
If $P:\C^N\rightarrow\R$ is a non-negative, real valued polynomial of
order $d$, such that the homogeneous terms of order $d$ are strictly
positive outside the origin, then $P$ has a minimum.
\end{lemma}
From lemmas \ref{lemmapol1} and \ref{lemmapol2} one gets:
\begin{theorem}\cite{CGR2}
Every coadjoint orbit of a compact semi-simple Lie group $G$
admits a preferred invariant symplectic connection.
\end{theorem}
In the case of a coadjoint orbit of $SU(3)$, one can prove by
direct calculation that this preferred connection is unique.

\section{Ricci-type connections} \label{Ricci type}
 
Ricci-type connections were defined in Definition \ref{def:Riccitype};
they are symplectic connections whose curvature tensor is entirely
determined by the Ricci tensor, i.e.\ for which $W^\nabla=0$.

\subsection{Some properties of the curvature of a Ricci-type connection}

Let $(M,\omega)$ be a smooth symplectic manifold of dim $2n$ ($n
\geq 2$) and let $\nabla$ be a smooth Ricci-type symplectic connection. 
The following results follow directly from the definition 
(and Bianchi's second identity).
 \begin{lemma}\label{properties}
\cite{CGR}
The curvature endomorphism has the form
\begin{equation}\label{curv}
R^\nabla(X,Y)=-\frac{1}{2(n+1)}[-2\omega (X,Y)\rho^\nabla-\rho^\nabla Y
\otimes \underline{X} + \rho^\nabla X \otimes \underline{Y} - X \otimes
\underline{\rho^\nabla Y} + Y \otimes \underline{\rho^\nabla X}]
\end{equation}
 where $\underline{X}$ denotes the 1-form $i(X) \omega$ (for $X$ a
 vector field on $M$) and where, as before, $\rho^\nabla$ is the
 endomorphism associated to the Ricci tensor by $ r^\nabla(U,V) =\omega
 (U, \rho^\nabla V)$.

\noindent Furthermore: 
\begin{itemize} 
\item[(i)] there exists a vector field $U^\nabla$ such that
\begin{equation}\label{u}
\nabla_X \rho^\nabla = -\dfrac{1}{2n+1} [X \otimes \underline{U}^\nabla
+ U^\nabla \otimes \underline{X}];
\end{equation}
\item[(ii)] there exists a function $f^\nabla$ such that
\begin{equation}\label{f}
\nabla_X U^\nabla = -\dfrac{2n+1}{2(n+1)} (\rho^\nabla)^2 X + f^\nabla X;
\end{equation}
\item[(iii)] there exists a real number $K^\nabla$ such that
\begin{equation}\label{k}
tr (\rho^\nabla)^2 + \dfrac{4(n+1)}{2n+1}f^\nabla = K^\nabla.
\end{equation}
\end{itemize}
\end{lemma}
\begin{corollary}
Any Ricci-type symplectic connection is preferred
\end{corollary}
The fact that $\cyclic_{XYZ}\nabla_Xr^\nabla(Y,Z)=0$
follows immediately from (\ref{u}).

\subsection{A construction by reduction}\label{section:reduction}

Consider the  manifold $M=\R^{2n+2}$ with its standard symplectic
structure $\Omega'$.\\
Let $A$ be a non-zero element in the symplectic Lie algebra
$sp(\R^{2n+2},\Omega')$.\\
Let $\Sigma_A$ be
the closed hypersurface $\Sigma_A \subset {\R}^{2n+2}$ defined by
$$
\Sigma _A=\{ x\in \R^{2n+2}\vert \Omega'(x,Ax)=1\}.
$$ 
(In order for $\Sigma_A$ to be non-empty we replace, if
necessary, $A$, by $-A$.)

The $1$-parameter subgroup $\exp tA$ of the symplectic group acts on
$\R^{2n+2}$, preserving $\Omega'$ and $\Sigma_A$; [the corresponding
fundamental vector field $A^*$ on $\R^{2n+2}$ (defined by
$A^*_x:=\frac{d}{dt}\exp -tA x_{\vert_0}=-Ax$) is Hamiltonian, i.e.
$i(A^*)\omega =dH_A $, with $H_A(x)=\frac{1}{2}\Omega(x,Ax)$ and
$\Sigma_A$ is a level set of this Hamiltonian].

We shall consider the {\bf reduced space} $ M^{red}:=\Sigma_A /\{\exp tA
\suchthat t \in \R\} $ with the canonical projection
$\pi:\Sigma_A\rightarrow M^{red}$. This can always be locally defined as
follows.

\noindent Since the vector field $Ax$ is nowhere 0 on $\Sigma_A$, for
any $x_0 \in \Sigma_A$, there exist:

-- a neighbourhood $U_{x_0}(\subset\Sigma_A)$, 

-- a ball $D^{red}\subset \R^{2n}$ of radius $r_0$, centred at the origin,

-- a real interval $I=(-\epsilon,\epsilon)$ 

-- and a diffeomorphism 
\begin{equation}
\chi:D^{red}\times I \to U_{x_0}
\end{equation}
such that $\chi(0,0)=x_0$ and $\chi(y,t) =\exp -tA(\chi(y,0))$. We
shall denote
$$
\pi: U_{x_0} \to D^{red} \quad  \pi = p_1
\otimes \chi^{-1}.
$$
The space $D^{red}$ is a local version of the Marsden--Weinstein
reduction of $\Sigma_A$ around the point $x_0$.

If $x \in \Sigma_A$, $T_x\Sigma_A =\ \rangle Ax \langle^{\perp}$, where
$\rangle v_1,\dots,v_p \langle$ denotes the subspace spanned by
$v_1,\dots,v_p$ and $~^\perp$ denotes the orthogonal relative to
$\Omega'$; let ${\mathcal H}_x (\subset T_x \Sigma_A) = \rangle x, Ax
\langle^\perp$; then ${\pi_*}_x$ defines an isomorphism between
${\mathcal H}_x$ and the tangent space $T_yD^{red}$ for $y=\pi(x)$.

A {\bf{reduced symplectic form}} on $D^{red}$, $\omega^{red}$, is
defined by
\begin{equation}
\omega^{red}_{y=\pi(x)}(X,Y):=\Omega'_x({\overline{X}}_x,{\overline{Y}}_x)
\end{equation}
where $\overline{Z}$ denotes the horizontal lift of $Z\in T_yD^{red}$;
i.e. $\overline{Z}\in {\mathcal H}_x$ and ${\pi_*}_x(\overline{Z})=Z$.

Let ${\nabla}$ be the standard flat symplectic affine connection on
$\mathbb R^{2n+2}$. The {\bf{reduced symplectic connection}}
$\nabla^{red}$ on $D^{red}$ is defined by
\begin{equation}
(\nabla^{red}_X Y)_y:=\pi_{*_x}(\nabla_{\overline{X}}{\overline{Y}}
-\Omega'(A{\overline{X}},{\overline{Y}})x +\Omega'
({\overline{X}},{\overline{Y}})Ax)
\end{equation}

\begin{prop}\label{ured}
\cite{BaguisCahen} The manifold $(D^{red}, \omega^{red})$ is a
symplectic manifold and $\nabla^{red}$ is a symplectic connection of
Ricci-type on it.

Furthermore, a direct computation shows that the corresponding
$\rho^{\nabla^{red}}, U^{\nabla^{red}}$ and $f^{\nabla^{red}}$ are given
by:
\begin{eqnarray}
\overline{\rho^{\nabla^{red}} X}(x) &=& -2(n+1)\overline{A_x}{\bar{X}}\\[2mm]
\bar{U}^{\nabla^{red}}(x)&=&-2(n+1)(2n+1)\overline{A^2_x}x\\[2mm]
(\pi^*f^{\nabla^{red}})(x)&=&2(n+1)(2n+1)\Omega'(A^2x,Ax)
\end{eqnarray}
where $\overline{A^k_x}$ is the map induced by $A^k$ with values in
${\mathcal H}_x$:
$$
\overline{A^k_x}(X)= A^k X + \Omega'(A^k X, x) Ax - \Omega'(A^k X, Ax) x .
$$
\end{prop}

\subsection{Local models for Ricci-type connections}

We have seen that given a Ricci-type symplectic connection $\nabla$ on a
symplectic manifold $(M,\omega)$ the curvature $R^\nabla$ is entirely
determined by $\rho^\nabla$ (\ref{curv}); its covariant derivative
$\nabla R^\nabla$ is thus determined by $\nabla \rho^\nabla$ which in
turn is determined by the vector field $U^\nabla$ (\ref{u}). The second
covariant derivative $\nabla^2 R^\nabla$ is determined by $\nabla
U^\nabla$ hence by $\rho^\nabla$ and $f^\nabla$ (\ref{f}). Since
$f^\nabla$ satisfies equation (\ref{k}), all successive covariant
derivatives of the curvature tensor are determined by
$\rho^\nabla,~U^\nabla$ and $ K^\nabla$.

\begin{corollary}
Let $(M,\omega)$ be a smooth symplectic manifold of dimension $2n$ $(n
\geq 2)$ and let $\nabla$ be a smooth Ricci-type connection. Let $p_0
\in M$; then the curvature $R^\nabla_{p_0}$ and its covariant
derivatives $(\nabla^k R^\nabla)_{p_0}$ (for all $k$) are determined by
$(\rho^\nabla_{x_0}, U^\nabla_{x_0}, K^\nabla)$.
\end{corollary}

\begin{corollary}
Let $(M, \omega, \nabla)$ (resp. $(M', \omega', \nabla')$) be two
symplectic manifolds of the same dimension $2n$ $(n \geq 2)$ each of
them endowed with a symplectic connection of Ricci-type. Assume that
there exists a linear map $b: T_{x_0} M \to T_{x'_0}M'$ such that 
\begin{itemize}
\item[(i)] $b^* \omega'_{x'_0} = \omega_{x_0}$, 
\item[(ii)] $b u^\nabla_{x_0}=u^{\nabla'}_{x'_0}$, 
\item[(iii)] $b \circ \rho^\nabla_{x_0} \circ b^{-1} =
\rho^{\nabla'}_{x'_0}$.
\end{itemize}
Assume further that $K^\nabla = K^{\nabla'}$. Then the manifolds are
locally affinely symplectically isomorphic, i.e.\ there exists a
normal neighbourhood of $x_0$ (resp. $x'_0$) $U_{x_0}$ (resp.\
$U'_{x'_0}$) and a symplectic affine diffeomorphism $\varphi: (U_{x_0},
\omega, \nabla) \to (U'_{x'_0}, \omega', \nabla')$ such that
$\varphi(x_0) = x'_0$ and $\varphi_{*x_0} = b$.
\end{corollary}

\noindent In case the symplectic manifold and the connection are real 
analytic, this follows from classical results, see for instance Theorem 7.2
and Corollary 7.3 in Kobayashi-Nomizu Volume 1 \cite{KN}. However, 
connections of Ricci-type are {\em always} real analytic, as we shall
see in Section~\ref{section:SpecialSymplectic}.

\begin{theorem} \label{thm:Ricci model}
Any symplectic manifold with a Ricci-type connection is locally
symplectically affinely isomorphic to the symplectic manifold with a
Ricci-type connection obtained by a local reduction procedure around
$e_0= (1,0,\dots,0)$ from a constraint surface $\Sigma_A$ defined by a
second order polynomial $H_A$ for $A\in sp(\R^{2n+2},\Omega')$ in the
standard  symplectic manifold $(\R^{2n+2},\Omega')$ endowed with the
standard flat connection.
\end{theorem}
Indeed if $p\in M$ and if $ \xi$ is a symplectic frame at $p$
[i.e. $\xi: ({\R}^{2n},\Omega^{(2n)})\rightarrow (T_{p},\omega_{p})$  is a
symplectic isomorphism of vector spaces], one defines
\begin{equation}\label{eq:tildeu}
\tilde{u}(\xi)=(\xi)^{-1}~U^\nabla(p),
\quad \tilde{\rho}(\xi)=(\xi)^{-1}~\rho^\nabla(p)~\xi
\end{equation}
and
\begin{equation}\label{eq:tildeA}
{\tilde{A}}(\xi)= \left( \begin{array}{ccc} 0 &\dfrac{
f(p)}{2(n+1)(2n+1)}
&\dfrac{-{\underline{\tilde{u}(\xi)}}}{2(n+1)(2n+1)}\\[2mm]
1 &0 &0\\[2mm]
0 &\dfrac{- \tilde{u}(\xi)}{2(n+1)(2n+1)}
&\dfrac{-\tilde{\rho}(\xi)}{2(n+1)}
\end{array}
\right)
\end{equation}
with ${\underline{u}}:=\Omega'(u,\cdot)$ and one looks at the
reduction for this $A={\tilde{A}}(\xi)$.

\subsection{Global models for Ricci-type connections} \label{globalRicci}

\begin{theorem}\cite{CahenGuttSchwach}
If $(M,\omega, \nabla)$  is of Ricci type with $M$ simply connected
there exists $(P,\omega^P)$ symplectic of dimension $2$ higher with a
flat connection $\nabla^P$ so that $(M,\omega, \nabla)$ is obtained from
$(P,\omega^P,\nabla^P)$ by reduction.
\end{theorem}

The manifold $P$ is obtained as the product $P=N\times \R$ of a contact
manifold $N$ and the real line $R$. The manifold $N$ is the holonomy
bundle over $M$ corresponding to a connection defined on the
$Sp({\mathbb{R}}^{2n+2},{\Omega'})$-principal bundle
$$
B'(M)=
B(M)\times_{Sp({\mathbb{R}}^{2n},{\Omega})}Sp({\mathbb{R}}^{2n+2},{\Omega'})
$$
with projection $\pi':B(M)'\to M$, where $B(M)\stackrel{\pi}{\to} M$ 
is the $Sp(\R^{2n},\Omega)$-principal bundle of symplectic frames
over $M$ and where we inject the symplectic group
$Sp(\R^{2n},\Omega)$ into $Sp({\mathbb{R}}^{2n+2},{\Omega'})$~
as the set of matrices
$$
\tilde{j}(A)=\left( \begin{array}{cc} I_2 &0\\0 &A
\end{array} \right)\qquad A\in Sp(\R^{2n},\Omega).
$$
The connection $1$-form $\alpha'$ on $B'(M)$ is characterised by the fact that
$$
\alpha'_{[\xi,1]}([\overline{X}^{hor},0])=\alpha_{\xi}(\overline{X}^{hor}).
$$
where
\begin{equation}\label{ref:alpha}
 \quad \alpha_{\xi}(\overline{X}^{hor})  = \left(
\begin{array}{ccc} 
0 &\dfrac{-\omega_x(u,X)}{2(n+1)(2n+1)}
&\dfrac{-\underline{\widetilde{\rho(X)}(\xi)}}{2(n+1)}\\[2mm]
0 & 0 & -\underline{\tilde{X}(\xi)}\\[2mm]
\tilde{X}(\xi) &\dfrac{-\widetilde{\rho(X)}(\xi)}{2(n+1)} &0
\end{array}
\right)
\end{equation}
where $X\in T_xM$ with $x=\pi(\xi)$ and $\overline{X}^{hor}$ 
is the horizontal lift of $X$ in $T_{\xi}B(M)$).

The equations satisfied by a Ricci-type connection imply that
the curvature $2$-form of the connection $1$-form $\alpha'$ is equal to
$ -2\tilde{A'} {\pi'}^*\omega$ where $\tilde{A'}$ is the unique
$Sp({\mathbb{R}}^{2n+2},{\Omega'})$-equi\-var\-iant extension of $\tilde{A}$
to $B'(M)$; and this curvature $2$-form is invariant by parallel transport 
($d^{\alpha'}curv(\alpha')=0$).\\
Thus the holonomy algebra of $\alpha'$ is of dimension $1$.
Assume $M$ is simply connected. The holonomy bundle
of $\alpha'$ is a circle or a line bundle over $M$,
$N\stackrel{\pi'}{\to} M$. This bundle has a natural
contact structure $\nu$ given by the restriction to 
$N\subset B(M)'$ of the $1$-form $-\half \alpha'$ (viewed as
real valued since it is valued in a $1$-dimensional algebra).
One has $d\nu={\pi'}^*\omega$.

The symplectic manifold with connection $(P,\omega^P,\nabla^P)$
is then obtained by an induction procedure that we shall
expose in a more general setting in Section \ref{section:induction}.

\section{Special symplectic connections} \label{section:SpecialSymplectic}

The striking rigidity results from Section~\ref{Ricci type} on Ricci-type
connections turn out to be a special case of a much more general phenomenon. 
As we saw, a connection of Ricci-type can be obtained by a symplectic 
reduction of a symplectic vector space with a flat symplectic connection. This
implies, for example, that the local moduli space of such connections is finite
dimensional.

As it turns out, this is merely a special case of a much broader
phenomenon. Indeed, there are many more geometric structures which can be
characterised in similar terms. For this, we call a symplectic connection on 
$(M, \om)$ with $\dim M \geq 4$ {\em special symplectic} if it belongs to one 
of the following classes.

\bi
\item {\bf Connections of Ricci-type} (cf.\ Section~\ref{Ricci type})

\item {\bf Bochner--K\"ahler and Bochner-bi-Lagrangian connections}

If the symplectic form is the K\"ahler form of a (pseudo-)K\"ahler metric, 
then its curvature decomposes into the Ricci curvature and the Bochner 
curvature (\cite{Bo}). If the latter vanishes, then (the Levi-Civita 
connection of) this metric is called Bochner--K\"ahler. 

Similarly, if the manifold is equipped with a bi-Lagrangian structure, i.e. 
two complementary Lagrangian distributions, then the curvature of a 
symplectic connection for which both distributions are parallel decomposes 
into the Ricci curvature and the Bochner curvature. Such a connection is 
called Bochner-bi-Lagrangian if its Bochner curvature vanishes. 

For results on Bochner--K\"ahler and Bochner-bi-Lagrangian connections, see 
\cite{Bochner} and \cite{K} and the references cited therein.

\item {\bf Connections with special symplectic holonomy}

A symplectic connection is said to have {\em special symplectic
holonomy} if its holonomy is contained in a proper absolutely
irreducible subgroup of the symplectic group.

The special symplectic holonomies have been classified in \cite{MS} and
further investigated in \cite{Br1}, \cite{CMS}, \cite{Habil}, \cite{S},
\cite{Advances}.
\ei

At first, it may seem unmotivated to collect all these structures in one
definition, but we shall provide ample justification for doing so.
Indeed, there is a beautiful link between special symplectic connections
and parabolic contact geometry.

For this, consider a simple Lie group $\G$ with Lie algebra $\g$. We say
that $\g$ is $2$-gradable, if $\g$ contains the root space of a long
root. This is equivalent to saying that there is a decomposition as a
graded vector space

\be \label{eq:decompose-g}
\g = \g^{-2} \oplus \g^{-1} \oplus \g^0 \oplus \g^1 \oplus 
\g^2,\ \ \ \ 
\mbox{and}\ \ \ \ [\g^i, \g^j] \subset \g^{i+j},
\ee
with $\dim \g^{\pm 2} = 1$. Indeed, there is a (unique) element
$H_{\al_0} \in [\g^{-2}, \g^2] \subset \g^0$ such that $\g^i$ is the
eigenspace of $\ad(H_{\al_0})$ with eigenvalue $i = -2, \ldots, 2$, and
any non-zero element of $\g^{\pm 2}$ is a long root vector.

Denote by $\p:= \g^0 \oplus \g^1 \oplus \g^2 \leq \g$ and let $\eP
\subset \G$ be the corresponding connected Lie subgroup. It follows that
the homogeneous space $\cC:= \G/\eP$ carries a canonical $\G$-invariant
contact structure which is determined by the $\Ad_\eP$-invariant
distribution $\g^{-1} \mod\, \p \subset \g/\p \cong T\cC$. In fact, we may
regard $\cC$ as the projectivisation of the adjoint orbit of a maximal
root vector. That is, we view $\cC \subset \P^o(\g)$ where $\P^o(V)$
denotes the set of oriented lines through $0$ of a vector space $V$, so
that $\P^o(V)$ is the sphere in $V$.

Each $a \in \g$ induces an action field $a^*$ on $\cC$ with flow $\T_a:=
\exp(\R a) \subset \G$ which hence preserves the contact structure on
$\cC$. Let $\cC_a \subset \cC$ be the open subset on which $a^*$ is
transversal to the contact distribution. There is a unique contact form
$\al \in \Om^1(\cC_a)$ determined by the equations that $\al(a^*) \equiv
1$. That is, $a^*$ is a Reeb vector field of the contact form $\al$.

We can cover $\cC_a$ by open sets $U$ such that the local quotient
$M_U:= \T_a^{loc} \backslash U$, i.e. the quotient of $U$ by a
sufficiently small neighbourhood of the identity in $\T_a$, is a
manifold. Then $M_U$ inherits a canonical symplectic structure $\om \in
\Om^2(M_U)$ such that $\pi^*(\om) = d\al$ for the canonical projection
$\pi: U \ra M_U$.

It is now our aim to construct a connection on $M_U$ which is
`naturally' associated to the given structure. For this, we let $\G_0
\subset \G$ be the connected subgroup with Lie algebra $\g^0 \leq \g$.
Since $\g^0 \leq \p$ and hence $\G_0 \subset \eP$, it follows that we
have a fibration
\be \label{eq:fibration}
\eP/\G_0 \longrightarrow \G/\G_0 \longrightarrow \cC = \G/\eP.
\ee

In fact, we may interpret $\G/\G_0:= \{(\al, v) \in T_p^*\cC \times
T_p\cC \mid p \in \cC, \al({\cal D}_p) = 0, \al(v) = 1\}$, where ${\cal
D} \subset T\cC$ denotes the contact distribution. Thus, given $a \in
\g$, then for each $p \in \cC_a$ we may regard the pair $(\al_p, a^*_p)$
from above as a point in $\G/\G_0$, i.e., we have a canonical embedding
$\imath: \cC_a \hookrightarrow \G/\G_0$.

Let $ \Gamma_a:= \pi^{-1}(\imath(\cC_a)) \subset \G$ where $\pi: \G
\ra \G/\G_0$ is the canonical projection. Then the restriction $\pi:
 \Gamma_a \ra \imath(\cC_a) \cong \cC_a$ becomes a principal
$\G_0$-bundle.

Consider the {\em Maurer-Cartan form} $\mu:= g^{-1} dg \in \Om^1(\G) \ot \g$
which we decompose according to (\ref{eq:decompose-g}) as $\mu = \sum_{i=-2}^2
\mu_i$ with $\mu_i \in \Om^1(\G) \ot \g^i$. Then we can show the following.

\begin{prop} \cite{CS}
Let $a \in \g$ be such that $\cC_a \subset \cC$ is non-empty, define the 
action field $a^* \in {\frak X}(\cC)$ and the principal $\G_0$-bundle $\pi: 
 \Gamma_a \ra \cC_a$ with $ \Gamma_a \subset \G$ from above. Then we 
have the following.
\bi
\item
The restriction of the components $\mu_0 + \mu_{-1} + \mu_{-2}$ of the 
Maurer-Cartan form to $\Gamma_a$ yields a pointwise linear isomorphism 
$T\Gamma_a \ra \g^0 \oplus \g^{-1} \oplus \g^{-2}$.
\item
There is a linear map $R: \g^0 \ra \Lambda^2 (\g^1)^* \ot \g^0$ and a smooth 
function $\rho: \Gamma_a \ra \g^0$ with the following property.
If we define the differential forms $\kappa \in \Om^1(\Gamma_a)$, $\th \in 
\Om^1(\Gamma_a) \ot \g^1$ and $\eta \in \Om^1(\Gamma_a) \ot \g^0$ by the 
equation
\[
\mu_0 + \mu_{-1} + \mu_{-2} = -2 \kappa\ \left(\frac12 e_{-2} + \rho\right) + 
\th + \eta
\]
for a fixed element $0 \neq e_{-2} \in \g^{-2}$, then the following equations
hold:
\be \label{eq:structurekappa}
d\kappa = \frac12 < e_{-2}, [\th, \th]>,
\ee
and
\be \label{eq:structureGamma} \ba{lll}
d\th + \eta \w \th & = & 0,\\ \\
d\eta + \frac12 [\eta, \eta] & = & R_{\rho} (\th \w \th).
\ea \ee
\ei
\end{prop}

Since the Maurer-Cartan form and hence $\kappa$, $\th$ and $\eta$ are
invariant under the left action of the subgroup $\T_a \subset \G$, we
immediately get the following

\begin{corollary} \cite{CS} 
On $\T_a \backslash \Gamma_a$, there is a coframing $\eta + \th \in
\Om^1(\T_a \backslash \Gamma_a) \ot (\g^0 \oplus \g^1)$ satisfying the
structure equations (\ref{eq:structureGamma}) for a suitable function
$\rho: \T_a \backslash \Gamma_a \ra \g^0$.
\end{corollary}

Thus, we could, in principle, regard $\th$ and $\eta$ as the tautological
and the connection $1$-form, respectively, of a connection on the
principal bundle $\T_a \backslash \Gamma_a \ra \T_a \backslash
\Gamma_a/\G_0$ whose curvature is represented by $R_\rho$. However,
$\T_a \backslash \Gamma_a/\G_0 \cong \T_a \backslash \cC_a$ will in
general be neither Hausdorff nor locally Euclidean, so the notion of a
principal bundle cannot be defined globally.

The way out of this difficulty is to consider {\em local} quotients
only, i.e., we restrict to sufficiently small open subsets $U \subset
\cC_a$ for which the local quotient $T_a^{loc} \backslash U$ {\em is} a
manifold. Clearly, $\cC_a$ can be covered by such open cells.

Moreover, if we describe explicitly the curvature endomorphisms $R_\rho$
for $\rho \in \g^0$, then one can show that -- depending on the choice of
the $2$-gradable simple Lie algebra $\g$ -- the connections constructed
above satisfy one of the conditions for a special symplectic connection
mentioned before.

More precisely, we have the following

\begin{theorem} \label{thmA} \cite{CS}
Let $\g$ be a simple $2$-gradable Lie algebra with $\dim \g \geq 14$,
and let $\cC \subset \P^o(\g)$ be the projectivisation of the adjoint
orbit of a maximal root vector. Let $a \in \g$ be such that $\cC_a
\subset \cC$ is non-empty, and let $\T_a = \exp(\R a) \subset \G$. If
for an open subset $U \subset \cC_a$ the local quotient $M_U =
\T_a^{loc} \backslash U$ is a manifold, then $M_U$ carries a special
symplectic connection.
\end{theorem}

The dimension restriction on $\g$ guarantees that $\dim M_U \geq 4$ and
rules out the Lie algebras of type $A_1$, $A_2$ and $B_2$.

The type of special symplectic connection on $M_U$ is determined by the
Lie algebra $\g$. In fact, there is a one-to-one correspondence between
the various conditions for special symplectic connections and simple
$2$-gradable Lie algebras. More specifically, if the Lie algebra $\g$ is
of type $A_n$, then the connections in Theorem~\ref{thmA} are
Bochner--K\"ahler of signature $(p,q)$ if $\g = \su(p+1,q+1)$ or
Bochner-bi-Lagrangian if $\g = \sl(n,\R)$; if $\g$ is of type $C_n$,
then $\g = \sp(n,\R)$ and these connections are of Ricci-type; if $\g$
is a $2$-gradable Lie algebra of one of the remaining types, then the
holonomy of $M_U$ is contained in one of the special symplectic holonomy
groups. Also, for two elements $a, a' \in \g$ for which $\cC_a, \cC_{a'}
\subset \cC$ are non-empty, the corresponding connections from
Theorem~\ref{thmA} are equivalent if and only if $a'$ is $\G$-conjugate
to a positive multiple of $a$.

If $\T_a \cong S^1$ then $\T_a \backslash \cC_a$ is an orbifold which
carries a special symplectic orbifold connection by Theorem~\ref{thmA}.
Hence it may be viewed as the ``standard orbifold model'' for (the
adjoint orbit of) $a \in \g$. For example, in the case of positive
definite Bochner--K\"ahler metrics, we have $\cC \cong S^{2n+1}$, and for
connections of Ricci-type, we have $\cC \cong \R\P^{2n+1}$. Thus, in
both cases the orbifolds $\T_a \backslash \cC$ are weighted projective
spaces if $\T_a \cong S^1$, hence the standard orbifold models $\T_a
\backslash \cC_a \subset \T_a \backslash \cC$ are open subsets of
weighted projective spaces.

Surprisingly, the connections from Theorem~\ref{thmA} exhaust {\em all}
special symplectic connections, at least locally. Namely we have the
following

\begin{theorem} \label{thmB} \cite{CS}
Let $(M, \om)$ be a symplectic manifold with a special 
symplectic connection of class $C^4$, and let $\g$ be the Lie algebra 
associated to the special symplectic condition as above.
\bi
\item
Then there is a principal $\hat \T$-bundle $\hat M \ra M$, where $\hat
\T$ is a one dimensional Lie group which is not necessarily connected,
and this bundle carries a principal connection with curvature $\om$.
\item
Let $\T \subset \hat \T$ be the identity component.
Then there is an $a \in \g$ such that $\T \cong \T_a \subset \G$, and a
$\T_a$-equivariant local diffeomorphism $\hat \imath: \hat M \ra \cC_a$ which 
for each sufficiently small open subset $V \subset \hat M$ induces a 
connection preserving diffeomorphism $\imath: \T^{loc} \backslash V \ra 
\T_a^{loc} \backslash U = M_U$, where $U:= \hat \imath(V) \subset \cC_a$ 
and $M_U$ carries the connection from Theorem~\ref{thmA}.
\ei
\end{theorem}

The situation in Theorem~\ref{thmB} can be illustrated by the following 
commutative diagram, where the vertical maps are quotients by the indicated 
Lie groups, and $\T \backslash \hat M \ra M$ is a regular covering.

\be \label{diagram2}
\xymatrix{
& \hat M \ar[d]^{\mbox{\footnotesize $\T$}} \ar[r]^{\hat \imath}
\ar[dl]_{\mbox{\footnotesize $\hat \T$}} & 
\cC_a \ar[d]^{\mbox{\footnotesize $\T_a$}} 
\\
M & \T \backslash \hat M \ar[r]^{\imath} \ar[l] & \T_a \backslash \cC_a 
}
\ee

In fact, one might be tempted to summarize Theorems~\ref{thmA} and
\ref{thmB} by saying that for each $a \in \g$, the quotient $\T_a
\backslash \cC_a$ carries a canonical special symplectic connection, and
the map $\imath: \T \backslash \hat M \ra \T_a \backslash \cC_a$ is a
connection preserving local diffeomorphism. If $\T_a \backslash \cC_a$
is a manifold or an orbifold, then this is indeed correct. In general,
however, $\T_a \backslash \cC_a$ may be neither Hausdorff nor locally
Euclidean, hence one has to formulate these results more carefully.

As consequences, we obtain the following

\begin{corollary} \label{corC}
All special symplectic connections of $C^4$-regularity are analytic, and
the local moduli space of these connections is finite dimensional, in
the sense that the germ of the connection at one point up to $3$rd order
determines the connection entirely. In fact, the generic special
symplectic connection associated to the Lie algebra $\g$ depends on
$(\rk(\g) - 1)$ parameters.

Moreover, the  Lie algebra ${\frak s}$ of {\em affine vector fields}, 
i.e., vector fields on $M$ whose flow preserves the connection, is 
 isomorphic to $stab(a)/(\R a)$ with $a \in
\g$ from Theorem~\ref{thmB}, where $stab(a) = \{ x \in \g \mid [x,a] =
0\}$. In particular, $\dim \s \geq \rk(\g) - 1$ with equality implying
that $\s$ is abelian.
\end{corollary}

When counting the parameters in the above corollary, we regard
homothetic special symplectic connections as equal, i.e. $(M, \om,
\nabla)$ is considered equivalent to $(M, e^{t_0} \om, \nabla)$ for all
$t_0 \in \R$.

We can generalize Theorem~\ref{thmB} and Corollary~\ref{corC} easily to
orbifolds. Indeed, if $M$ is an orbifold with a special symplectic
connection, then we can write $M = \hat \T \backslash \hat M$ where
$\hat M$ is a manifold and $\hat \T$ is a one dimensional Lie group
acting properly and locally freely on $\hat M$, and there is a local
diffeomorphism $\hat \imath: \hat M \ra \cC_a$ with the properties
stated in Theorem~\ref{thmB}.

There is a remarkable similarity between the cones $\cC_i \subset \g_i$,
$i = 1,2$, for the simple Lie algebras $\g_1:= \su(n+1,1)$ and $\g_2:=
\sp(n,\R)$. Namely, $\cC_1 = S^{2n+1}$ with the standard $CR$-structure,
and $\g_1$ is the Lie algebra of the group $\SU(n+1,1)$ of
$CR$-isomorphisms of $S^{2n+1}$ \cite{K}. On the other hand, $\cC_2 =
\R\P^{2n+1}$, regarded as the lines in $\R^{2n+2}$ with the
projectivised action of $\sp(n+1,\R)$ on $\R^{2n+2}$. Thus, $\cC_1$ is
the universal cover of $\cC_2$, so that the local quotients $\T_a
\backslash \cC_a$ are related. In fact, we have the following result.

\begin{prop} \cite{PS}
Consider the action of the $2$-gradable Lie algebras $\g_1:= \su(n+1,1)$
and $\g_2:= \sp(n+1,\R)$ on the projectivised orbits $\cC_1$ and
$\cC_2$, respectively. Then the following are equivalent.
\bi
\item
For $a_i \in \g_i$ the actions of $\T_{a_i} \subset \G_i$ on $\cC_i$ are
conjugate for $i = 1,2$,
\item
$a_i \in \uu(n+1)$ where $\uu(n+1) \subset \g_i$ for $i = 1,2$ via the two 
standard embeddings.
\ei
\end{prop}

This together with the preceding results yields the following

\begin{theorem} \cite{PS}
\bi
\item
Let $(M, \om, \nabla)$ be a symplectic manifold with a connection of
Ricci type, and suppose that the corresponding element $A \in
\sp(n+1,\R)$ from Theorem~\ref{thm:Ricci model} is conjugate to an
element of $\uu(n+1) \subset \sp(n+1,\R)$. Then $M$ carries a canonical
Bochner--K\"ahler metric whose K\"ahler form is given by $\om$.
\item
Conversely, let $(M, \om, J)$ be a Bochner-K\"ahler metric such that the
element $a \in \su(n+1,1)$ from Theorem~\ref{thmB} is conjugate to an
element of $\uu(n+1) \subset \su(n+1,1)$. Then $(M, \om)$ carries a
canonical connection of Ricci-type.
\ei
\end{theorem}

Note that in \cite{Bochner}, Bochner--K\"ahler metrics have been locally
classified. In this terminology, the Bochner--K\"ahler metrics in the
above theorem are called {\em Bochner--K\"ahler metrics of type I}.

\section{Symplectic twistor space and Ricci-type connections}

In this section we present a result of Vaisman \cite{Vaisman2} which
shows how the Ricci-type condition on the curvature of a symplectic
connection can be seen as an integrability condition for an associated
almost complex structure on the total space of a bundle over the
symplectic manifold. It is far from clear what is the geometrical
significance of this complex structure. Preliminary studies of its
properties have been made in the PhD theses of Albuquerque
\cite{Albuquerque} and Sti\'enon \cite{Stienon}.

\subsection{Compatible almost complex structures}

Let $(V,\Omega)$ be a finite-dimensional real symplectic vector space.
We denote by $Sp(V,\Omega)$ the real symplectic group of linear
transformations $g$ of $V$ which preserve $\Omega$. Any two symplectic
vector spaces of the same dimension are isomorphic, and $Sp(V,\Omega)$
acts freely and transitively on the set of isomorphisms from
$(V,\Omega)$ to any other symplectic vector space of the same dimension
by composition on the right.

The Lie algebra $\sp(V,\Omega)$ of $Sp(V,\Omega)$ consists of all
linear endomorphisms $\xi$ of $V$ satisfying
\[
\Omega(\xi u,v) + \Omega(u,\xi v)=0,\qquad \forall u,v \in V.
\]
This condition is equivalent to saying that $B_\xi(u,v)=\Omega( u,\xi v)$
defines a symmetric bilinear form $B_\xi$ on $V$. Conversely, any
symmetric bilinear form on $V$ defines an element of $\sp(V,\Omega)$. So
as vector spaces $\sp(V,\Omega)$ and $S^2V^*$ are isomorphic. In fact,
if we consider the natural actions of $Sp(V,\Omega)$ on both spaces we
have
\[
(g.B_\xi)(u,v) = B_\xi(g^{-1}u,g^{-1}v) = \Omega(g^{-1}u, \xi g^{-1}v) = 
\Omega(u,g\xi g^{-1}v) = B_{Ad_g\xi}(u,v)
\]
so that the adjoint representation of $Sp(V,\Omega)$ is isomorphic
to the natural representation of $Sp(V,\Omega)$ on $S^2V^* \cong S^2V$.

A \textit{compatible almost complex structure} $j$ on $(V,\Omega)$ is an
element of $\sp(V,\Omega)$ with $j^2 = -1$ such that $B_j$ is positive
definite and we denote the set of compatible almost complex structures
by $j(V,\Omega)$. $Sp(V,\Omega)$ acts transitively on $j(V,\Omega)$ by
conjugation. If we drop the positivity condition, then $Sp(V,\Omega)$
has a finite number of orbits distinguished by the signature of $B_j$.
For our purpose any orbit would do, we choose the positive orbit for
convenience.

If $j$ is a compatible almost complex structure on $(V,\Omega)$ then $V$
becomes a complex vector space by defining $(a+ib)v = av + bjv, \forall
a,b \in \R$. Further $\langle u,v \rangle_j = B_j(u,v) - i \Omega(u,v)$
defines a positive Hermitean structure on $V$ as a complex vector space.
Under the action of $Sp(V,\Omega)$ above, it is clear that the
stabiliser of $j$ is the unitary group $U(V,\Omega,j)$ of this Hermitean
structure so that $j(V,\Omega) \cong Sp(V,\Omega)/U(V,\Omega,j)$ as 
homogeneous spaces.

\subsection{Geometry of $j(V,\Omega)$}

Since $j(V,\Omega)$ is a homogeneous space, the tangent space at $j$ can
be identified with the quotient of Lie algebras
$\sp(V,\Omega)/\jhru(V,\Omega,j)$. $\jhru(V,\Omega,j)$ consists of all
elements of $\sp(V,\Omega)$ which commute with $j$ so the set of
elements which anticommute with $j$ form a complementary subspace
\[
\jhrm_j = \{\xi \in \sp(V,\Omega) \suchthat \xi j + j \xi = 0\}.
\]
Since $j\xi j \in \sp(V,\Omega)$ for $\xi  \in \sp(V,\Omega)$ we have
a decomposition
\[
\xi = \frac 12 (\xi - j \xi j) + \frac 12 (\xi + j \xi j) 
\]
with $\xi - j \xi j \in \jhru(V,\Omega,j)$ and $\xi + j \xi j \in \jhrm_j$.
This shows that
$\sp(V,\Omega) = \jhru(V,\Omega,j) +\jhrm_j$
is a direct sum and so we have a natural isomorphism
$T_j j(V,\Omega) \cong \m_j$.

The product of two endomorphisms which anticommute with $j$ must commute
with $j$ so $[\jhrm_j,\jhrm_j] \subset \jhru(V,\Omega,j)$ and hence
$Sp(V,\Omega)/U(V,\Omega,j)$ is a symmetric space. In fact this is a
realisation of the Siegel domain for the symplectic group. The canonical
connection is torsion free and any invariant object is parallel.

If we take $\xi \in \jhrm_j$ then $j\xi = \frac12 [j,\xi] \in \jhrm_j$ and, as
$j$ varies, $\xi \mapsto j\xi$ defines an almost complex structure on
$Sp(V,\Omega)/U(V,\Omega,j)$ which is clearly $Sp(V,\Omega)$-invariant.
Being invariant it is parallel and being parallel for a torsion-free
connection it is integrable and gives the invariant complex structure of
$j(V,\Omega)$ as a Hermitean symmetric  space.

If $j$ is an endomorphism of $V$ with $j^2=-1$ then $j$ extends complex
linearly to the complexification $V^{\C}$ and has $\pm i$ as eigenvalues
with eigenspaces $V^\pm$. Then $V^{\C} = V^+ + V^-$ is a direct sum and
the projections $j^{\pm}$ onto these subspaces are given by $j^\pm =
\frac12(1\mp ij)$.

Recall that $\sp(V,\Omega)$ coincides with $S^2V^* \cong S^2V$ as
representations. $j$ acts on both of these as an element of the Lie
algebra, and it is clear that the eigenvalues of $j$ on $S^2V^{\C}$ will
be sums of eigenvalues of $j$ on $V^{\C}$. Thus they are $0, \pm2i$ and
the same must hold on $\sp(V,\Omega)^{\C}$. The zero eigenspace
corresponds with $\jhru(V,\Omega,j)^{\C}$ and the $\pm2i$ eigenspaces split
$\jhrm_j^{\C}$ into $\jhrm_j^+$ and $\jhrm_j^-$. These are the $\pm i$
eigenspaces of multiplication by $j$ on the left which defined the
complex structure. Thus $\jhrm_j^+$ is the space of $(1,0)$ tangent vectors
to $j(V,\Omega)$ at $j$. Since there is no $4i$ eigenspace on
$\sp(V,\Omega)^{\C}$ we have $[\jhrm_j^+,\jhrm_j^+]=0$ and since there is no
$3i$ eigenspace on $V^{\C}$, $\jhrm_j^+ V^+ = 0$. Note that $\jhrm_j^- V^+ =
V^-$ so the condition $\jhrm_j^+ V^+ = 0$ is a compatibility condition
between the complex structures on $V$ and $j(V,\Omega)$ which
distinguishes the two possible invariant complex structures on
$j(V,\Omega)$ and allows us to choose one in preference to the other.

\subsection{The bundle $J(M,\omega)$ of almost complex structures}

We fix a symplectic vector space $(V,\Omega)$ and consider a symplectic
manifold $(M,\omega)$ of the same dimension as $V$. A \textit{symplectic
frame} at $x \in M$ is a symplectic isomorphism $p$ from $(V,\Omega)$ to
$(T_xM,\omega_x)$. We denote the set of symplectic frames at $x$ by
$Sp(M,\omega)_x$ and the disjoint union of these over $x \in M$ by
$Sp(M,\omega)$. With the obvious projection map $\pi \colon Sp(M,\omega)
\to M$, $Sp(M,\omega)$ is the \textit{symplectic frame bundle} of
$(M,\omega)$. It is a principal $Sp(V,\Omega)$ bundle with
$Sp(V,\Omega)$ acting on the right by composition.

Denote by $J(M,\omega)$ the bundle on $M$ whose fibre at $x$ is
$j(T_xM,\omega_x)$. Smooth sections of $J(M,\omega)$ are almost complex
structures on $M$ compatible with the symplectic structure so we call
$J(M,\omega)$ the \textit{bundle of compatible almost complex structures}.
$J(M,\omega)$ can also be viewed as the associated bundle
$Sp(M,\omega)\times_{Sp(V,\Omega)}j(V,\Omega)$ with fibre $j(V,\Omega)$.
Since the latter is homogeneous, $J(M,\omega)$ can be identified with
the quotient space $Sp(M,\omega)/U(V,\Omega,j)$ with the identification
given by
\begin{equation}\label{SympFB:quot}
Sp(M,\omega)/U(V,\Omega,j) \ni p.U(V,\Omega,j) \longleftrightarrow
p\circ j \circ p^{-1} \in J(M,\omega).
\end{equation}

$Sp(M,\omega)$ has a $V$-valued $1$-form $\theta$ called the 
\textit{soldering form} defined by
\[
\theta_p(X) = p^{-1} (\pi_* X)
\]
and a \textit{vertical bundle} $\jhrV$ with
\[
\jhrV_p = \Ker \pi_* \colon T_p Sp(M,\omega) \to T_{\pi(p)}M
\]
so $\jhrV$ is the kernel of $\theta$. A $1$-form on $Sp(M,\omega)$ is
called \textit{horizontal} if it vanishes on $\jhrV$\/. The components
of $\theta$ span the horizontal forms pointwise.

If $\xi \in \sp(V,\Omega)$ then $\widetilde{\xi}$ is the \textit{vertical
vectorfield} on $Sp(M,\omega)$ defined by:
\[
\widetilde{\xi}_p = \left.\frac{d}{dt} p \circ \exp t\xi \right|_{t=0}.
\]
$(p,\xi) \mapsto \widetilde{\xi}_p$ identifies the trivial bundle
$Sp(M,\omega)\times\sp(V,\Omega)$ with $\jhrV$. Note that the pull-back
bundle $\pi^{-1}TM$ to $Sp(M,\omega)$ is also trivial via the map $(p,X)
\mapsto (p, p^{-1}X)$ so $\theta$ gives a map $TSp(M,\omega) \to
Sp(M,\omega)\times V$ and an exact sequence of bundles on $Sp(M,\omega)$
\[
 0 \map{} Sp(M,\omega)\times\sp(V,\Omega) \map{}
 TSp(M,\omega) \map{} Sp(M,\omega)\times V \map{} 0.
\]
Splitting this exact sequence trivialises the tangent bundle of
$Sp(M,\omega)$ and makes calculations with forms particularly easy.

\subsection{Almost symplectic connections}

A principal connection in $Sp(M,\omega)$ is an \textit{almost symplectic
connection} $\nabla$ in $TM$ (a symplectic connection additionally has
vanishing torsion). It is given by an $\sp(V,\Omega)$-valued $1$-form
$\alpha$ on $Sp(M,\omega)$ which satisfies $\alpha({\widetilde\xi}) =
\xi$ and $g^*\alpha = \Ad_g \alpha$ for all $\xi \in \sp(V,\Omega)$ and
$g \in Sp(V,\Omega)$. It follows that $\H = \Ker \alpha$ is a
complementary subbundle to $\jhrV$ which we call the \textit{horizontal
bundle} determined by $\alpha$. $\nabla$ and $\alpha$ are related as
follows: if $U$ is an open set in $M$ on which there is a local section
$s$ of $Sp(M,\omega)$ then $sv$ is a vectorfield on $U$ for any fixed
$v$ in $V$ and $s^*\alpha$ is an $\sp(V,\Omega)$-valued $1$-form on $U$
which are related by
\[
\nabla_X (sv) = s ((s^*\alpha)(X)v)
\]
or, more compactly,
\[
s^{-1} \nabla s = s^* \alpha
\]
with the appropriate interpretation of the two sides of the equation.

The components of the forms $\alpha,\theta$ are linearly independent and
span the cotangent spaces of $Sp(M,\omega)$ at each point so that
$Sp(M,\omega)$ is parallelisable. This makes the differential geometry
on $Sp(M,\omega)$ especially simple and is the reason why we shall do
calculations on $Sp(M,\omega)$ in preference to $M$ or $J(M,\omega)$.

The torsion and curvature of $\nabla$ lift to $Sp(M,\omega)$ as
horizontal $V$-{} and $\sp(V,\Omega)$-valued $2$-forms $\tau^\nabla$ and
$\rho^\nabla$ given
by
\[
\tau^\nabla_p(X,Y) = p^{-1}(T^\nabla_{\pi(p)}(\pi_*X,\pi_*Y)),\qquad
\rho^\nabla_p(X,Y) = p^{-1}\circ R^\nabla_{\pi(p)}(\pi_*X,\pi_*Y)\circ p.
\]
These forms can be computed from the corresponding principal connection
$\alpha$ by
\[
\tau^\nabla = d\theta + \alpha \wedge \theta, \qquad \rho^\nabla =
d\alpha + \frac12 [\alpha\wedge\alpha]
\]
where $[\alpha\wedge\alpha]$ denotes simultaneous wedge as forms and Lie
bracket of the values.

\subsection{Differential geometry on $J(M,\omega)$}

Let $\pi$ be the projection $J(M,\omega)\to M$, then $\pi$ is a
submersion so $d\pi$ is surjective and we have an exact sequence
of vector bundles on $J(M,\omega)$
\[
0 \longrightarrow \jhrV \longrightarrow TJ(M,\omega) 
\stackrel{d\pi}{\longrightarrow} \pi^{-1}TM \longrightarrow 0
\]
where $\jhrV$ now denotes the vertical bundle on $J(M,\omega)$. It follows
that if $j \in J(M,\omega)$ then $\jhrV_j$ consists of all elements $\xi$
of $\sp(T_xM,\omega_x)$ which anticommute with $j$ where $x=\pi(j)$. We
identify $\jhrV_j$ with this subspace via the map 
\[
\sp(T_xM,\omega_x) \ni \xi \longleftrightarrow
\left.\frac{d}{dt} \exp t \xi \circ j \circ \exp -t \xi\right|_{t=0}.
\]

$\End(\pi^{-1}TM)$ has a tautological section which we denote by $\jhrJ$
whose value at $j$ is $j \in \End (T_{\pi(j)}M)$. If we denote by
$\sp(M,\omega)$ the bundle of Lie algebras whose fibre at $x$ is
$\sp(T_xM,\omega_x)$ then $\jhrJ$ is a section of $\pi^{-1}\sp(M,\omega)$.
The kernel of $\ad\jhrJ$ is a bundle $\jhru(M,\omega,\jhrJ)$ of unitary Lie
algebras on $J(M,\omega)$. The range of $\ad\jhrJ$ is the elements of
$\pi^{-1}\sp(M,\omega)$ which anticommute with $\jhrJ$ and so coincides
with the tangent space to the fibre, namely the vertical bundle $\jhrV$.
Thus
\[
\pi^{-1}\sp(M,\omega) = \jhru(M,\omega,\jhrJ) + \jhrV.
\]

If $\nabla$ is the almost symplectic connection on $TM$ corresponding to
a principal connection $\alpha$ then the horizontal distribution is
$U(V,\Omega,\jhrJ)$-invariant and hence projects to a horizontal
distribution $\jhrH^\nabla$ on $J(M,\omega)$. Our next objective is to
identify $\jhrH^\nabla$ directly without going via the frame
bundle.

The $\sp(V,\Omega)$ connection form $\alpha$ in $Sp(M,\omega)$ pulls
back to $\pi^*\alpha$ in $\pi^{-1}Sp(M,\omega)$ and then pulls
back to the reduction as $\alpha$. It has values in $\sp(V,\Omega)$
rather than the Lie algebra $\jhru(V,\Omega,j_0)$ so it is a Cartan
connection from the point of view of this reduction. We can split
$\alpha$ relative to the decomposition
\[
\sp(V,\Omega) = \jhru(V,\Omega,j_0) \oplus \jhrm_{j_0}
\]
as
\[
\alpha = \alpha^u + \alpha^m.
\]
Obviously, $\alpha^u$ is a principal connection in the bundle
$Sp(M,\omega) \longrightarrow J(M,\omega)$. $\alpha^m$ vanishes
on vertical vectors for this bundle so is the lift of a 
$Sp(M,\Omega)\times_{U(V,\Omega,j_0)} \jhrm_{j_0}$-valued $1$-form.
Since this associated bundle is $\jhrV$ we clearly have a $\jhrV$-valued
$1$-form on $J(M,\omega)$.

Let $\jhrP$ be the bundle map $TJ(M,\omega) \longrightarrow
\pi^{-1}\sp(M,\omega)$ with image $\jhrV$ and kernel $\jhrH^\nabla$.
Then $\jhrP$ anticommutes with $\jhrJ$. Note that $\jhrP$ can be viewed
as a $\pi^{-1}\sp(M,\omega)$-valued $1$-form.

\begin{prop}
The horizontal lift of $\jhrP$ to $Sp(M,\omega)$ is $\alpha^m$.
\end{prop}

The almost symplectic connection $\nabla$ induces a connection in
$\sp(M,\omega)$ and hence its pull-back $\pi^{-1}\nabla$ induces one in
$\pi^{-1}\sp(M,\omega)$. We can thus take the covariant derivative
$\pi^{-1}\nabla \jhrJ$ which is also a $\pi^{-1}\sp(M,\omega)$-valued
$1$-form. Moreover $\jhrJ^2=-1$ so $\pi^{-1}\nabla\jhrJ \jhrJ + \jhrJ
\pi^{-1}\nabla\jhrJ = 0$ and hence $\pi^{-1}\nabla \jhrJ$ is
$\jhrV$-valued. These two $1$-forms are related by

\begin{prop}
\[
\pi^{-1}\nabla\jhrJ = [\jhrP,\jhrJ].
\]
\end{prop}

\begin{proof}
$\jhrJ$ lifts to
$\pi^{-1}Sp(M,\omega)$ as $\widehat{\jhrJ}(j,p) = p\circ j \circ p^{-1}$
and so on the reduction is given by $\widehat{\jhrJ}(\sigma(p)) = j_0$. So
it is constant. On the reduction the covariant derivative is then
$\sigma^*(d\widehat{\jhrJ} + [\widehat{\pi}^*\alpha,\widehat{\jhrJ}]) =
[\alpha,j_0] = [\alpha^m,j_0]$.
\end{proof}

Finally we can describe the horizontal distribution on $J(M,\omega)$ directly.

\begin{corollary}
\[
\jhrH^\nabla = \{ X \in TJ(M,\omega) \suchthat \pi^{-1}\nabla_X \jhrJ = 0 \}.
\]
\end{corollary}

\subsection{The almost complex structures on $J(M,\omega)$}

If $(M,\omega)$ is a symplectic manifold the fibres of $J(M,\omega)$ are
diffeomorphic to $j(V,\Omega)$ which has an invariant complex structure,
which then transfers to each fibre. The tangent bundle to the fibres
$\jhrV$ thus has an endomorphism which gives the almost complex structure
of each fibre. Under the identification of $\jhrV$ with endomorphisms of
$\pi^{-1}TM$ this agrees with left multiplication by $\jhrJ$. 

If $\nabla$ is an almost symplectic connection on $(M,\omega)$ then the
horizontal distribution $\jhrH^\nabla$ is isomorphic via $d\pi$ with
$\pi^{-1}TM$ and the latter has an endomorphism $\jhrJ$. Since
$TJ(M,\omega) = \jhrV \oplus \jhrH^\nabla$ there is a unique
endomorphism $J^\nabla$ of $TJ(M,\omega)$ which coincides with
multiplication by $\jhrJ$ on $\jhrV$, and which satisfies $d\pi {\circ}
J^\nabla = \jhrJ {\circ} d\pi$. Clearly $J^\nabla$ is an almost complex
structure on $J(M,\omega)$. Thus we have shown

\begin{prop}
If $\nabla$ is an almost symplectic connection on the symplectic 
manifold $(M,\omega)$ then there is a unique almost complex structure
$J^\nabla$ on $J(M,\omega)$ which satisfies
\begin{itemize}
\item $\jhrP {\circ} J^\nabla = \jhrJ {\circ} \jhrP$, or equivalently,
$\pi^{-1}\nabla_{J^\nabla X}\jhrJ = \jhrJ \pi^{-1}\nabla_{X}\jhrJ$ for
all $X$;
\item $d\pi {\circ} J^\nabla = \jhrJ {\circ}
d\pi$.
\end{itemize}
\end{prop}

It is our goal to prove the
following theorem.

\begin{theorem}\label{theorem:twistor}
Let $\nabla$ be a connection on the symplectic manifold $(M,\omega)$.

\begin{enumerate} 
\item \label{jhr-torsion}If $\nabla$ is almost symplectic and $J^\nabla$
is integrable then there is a symplectic connection $\nablabar$ with
$J^{\nablabar} = J^\nabla$.

\item \label{jhr-injective}If $\nabla,\nablabar$ are symplectic
connections then $J^{\nablabar} =J^\nabla$ implies $\nabla=\nablabar$.

\item \label{jhr-integrable}$J^\nabla$ is integrable if and only if
$\nabla$ is of Ricci-type.
\end{enumerate}
\end{theorem}

The first step is to convert the integrability question into a condition 
on the frame bundle. This is possible since integrability of an almost
complex structure is equivalent to its $(1,0)$ forms generating
a $d$-closed ideal. Since the map $Sp(M,\omega) \to J(M,\omega)$ is a
submersion, the pull-back of forms is an injective map, so any relation
amongst forms on $J(M,\omega)$ will hold if and only if it holds
amongst their pull-backs to $Sp(M,\omega)$. Thus we will have integrability
if and only if a pointwise basis for the pull-backs of $(1,0)$ forms
generates a $d$-closed ideal on $Sp(M,\omega)$.

\begin{lemma}
The pull-backs to $Sp(M,\omega)$ of $(1,0)$ forms for $J^\nabla$ are
spanned pointwise by $\alpha^+ = j^+\alpha j^-$ and $\theta^+ =
j^+\theta$ if $\alpha$ is the connection form for $\nabla$.
\end{lemma}

See \cite{OBR} for a proof. An immediate consequence is

\begin{lemma}
$J^\nabla$ is integrable if and only if $j^+\tau^\nabla$ and
$j^+\rho^\nabla j^-$ are in the ideal generated by $\alpha^+, \theta^+$.
\end{lemma}

\begin{proof}
\[
j^+\tau^\nabla - d\theta^+ = j^+\alpha \wedge \theta = j^+(\alpha(j^+ +
j^-) \wedge \theta) =  j^+\alpha\wedge \theta^+ + \alpha^+\wedge \theta
\]
and
\[
j^+\rho^\nabla j^- - d\alpha^+ = j^+\rho^\nabla - d\alpha j^- =
\frac12j^+[\alpha,\alpha]j^- =   [\alpha^u,\alpha^+].
\]
This gives the Lemma immediately.
\end{proof}

\medskip

\noindent \textbf{Proof} of \tref{theorem:twistor} 
(\ref{jhr-torsion}):\ \ If $\nabla$ is almost symplectic with torsion
$T^\nabla$ then it is easily verified that the connection $\nablabar$
defined by
\[
\omega(\nablabar_XY, Z) = \omega(\nabla_XY, Z) - \frac12
\omega(T^\nabla(X,Y),Z) + \frac16 \omega(X, T^\nabla(Y,Z)) + \frac16
\omega(Y, T^\nabla(X,Z))
\]
is a symplectic connection. The connection forms $\alpha$
and $\widebar{\alpha}$ on $Sp(M,\omega)$ are related by
\[
\Omega(\widebar{\alpha}u, v) = \Omega(\alpha u, v) - \frac12
\Omega(\tau^\nabla(\theta,u),v) + \frac16 \Omega(\theta, \tau^\nabla(u,v)) +
\frac16 \Omega(u, \tau^\nabla(\theta,v)).
\]
Hence
\begin{eqnarray*}
\Omega(\widebar{\alpha}^+u, v) &=& \Omega(\widebar{\alpha}j^-u,
j^-v)\\
 &=& \Omega(\alpha^+ u, v) - \frac12\Omega(j^+\tau^\nabla(\theta,j^-u),v)\\
&&\qquad\qquad + \frac16 \Omega(\theta,j^+\tau^\nabla(j^-u,j^-v))
 + \frac16 \Omega(u, j^+\tau^\nabla(\theta,j^-v)).
\end{eqnarray*}
Integrability implies that $j^+\tau^\nabla( j^-u, j^-v)=0$ from which
one sees that $ j^+\tau^\nabla(\theta, j^-v) = j^+\tau^\nabla(\theta^+,
j^-v)$ and hence that $\widebar{\alpha}^+$ is a combination of
$\alpha^+$ and $\theta^+$ so the $(1,0)$ forms of $J^\nabla$ and
$J^\nablabar$ agree and hence the almost complex structures are the
same.

\medskip

\noindent \textbf{Proof} of \tref{theorem:twistor}
(\ref{jhr-injective}):\ \ Suppose $\nabla$ and $\nablabar$ are both
symplectic connections defining the same almost complex structures
$J^\nabla = J^\nablabar$, which will be the case if and only if the
components of $\widebar{\alpha}^+$ are combinations of the components of
$\alpha^+$ and $\theta^+$. If $\nablabar = \nabla + B$ and we lift $B$
to $Sp(M,\omega)$ as
\[
\beta_p (X)  =  p^{-1}{\circ} B_{\pi(p)}(\pi_*X) {\circ} p
\]
then $\widebar{\alpha} = \alpha +\beta $ so $J^\nablabar = J^\nabla$ if
and only if $\beta^+$ is a combination of the components of $\alpha^+$
and $\theta^+$, and since $\beta$ is horizontal, this means $\beta^+$
vanishes on $(0,1)$ vectors. Hence $j^+ B_x(j^-X)j^-Y = 0$ for every $j
\in j(T_xM,\omega_x)$.

Solving this condition can be converted to a problem in representation
theory by observing first that it can be re-expressed in terms of the
completely symmetric 3-tensor $\underline{B}(X,Y,Z) = \omega(B(X)Y,Z)$
as $\underline{B}(j^-X,j^-Y,j^-Z) = 0$. Since the set of such $j$ is
$Sp(V,\Omega)$-invariant, this says that $\underline{B}$ takes its
values in an invariant subspace in which there is no $3i$ eigenvalue for
any $j$ as an element of the Lie algebra. But $S^3V$ is irreducible and
the $3i$ eigenvalue does occur (on $S^3 V^+$) hence the only invariant
subspace with no $3i$ eigenvalue is $\{0\}$. Thus $B=0$ and hence
$\nablabar=\nabla$.

\medskip

\noindent \textbf{Proof} of \tref{theorem:twistor}
(\ref{jhr-integrable}):\ \ We again write the integrability condition in
terms of representation theory as in the proof of injectivity. The point
is that $j^+\rho^\nabla j^-$ consists of horizontal forms so is in the
ideal generated by $(1,0)$ forms if and only if $j^+R_x^\nabla(j^- X,
j^- Y) j^- Z = 0$ for all $j \in j(T_xM,\omega_x)$. This is again an
invariant condition, so it means $R^\nabla$ has to take its values in
the largest invariant subspace of the space of curvature tensors where
there is no $4i$ eigenvalue of $j$ as an element of the Lie algebra. We
know that there are two irreducible components to curvature, and one is
the Ricci component which is isomorphic to symmetric 2-tensors, so only
has eigenvalues $0$ and $\pm2i$. The other irreducible component does
have $4i$ eigenvalues and so for integrability the curvature must be
Ricci-type.

\subsection{Generalisations}

The fibres of $J(M,\omega)$ are large which is the reason the
integrability condition kills off so much of the curvature. We can try
to find subbundles with smaller fibres which are complex submanifolds of
the fibres of $J(M,\omega)$ and then play the same game with this
smaller bundle. For example we might have a $G$-structure $P \to M$ for
some subgroup $G$ of $Sp(V,\Omega)$ and a complex orbit $G/H$ of $G$ on
$j(V,\Omega)$. Then $P \times_G G/H = P/H$ will have an almost complex
structure for each choice of principal $G$-connection and a similar
calculation can be performed to determine when it is integrable. If $G$
is small then the curvature and torsion may break into several
irreducible pieces more of which may survive the integrability
condition. On the other hand as $G$ gets smaller it is harder to have
such a $G$-structure

Alternatively, for a given symplectic connection $\nabla$ we can look at
the zero-set of the Nijenhuis tensor of $J^\nabla$ on $J(M,\omega)$ for
a given symplectic connection $\nabla$. If a component of this set turns
out to be a complex manifold then it can be seen as a twistorial space
over $M$.

Both these approaches lead to interesting examples of twistor spaces in
the Riemannian case. Little is yet known in the symplectic case.

\subsection{The bundle $J(M, \om)$ for spaces of Ricci-type}%
\label{subsect:twistorreduction}

In view of \tref{theorem:twistor} it is interesting to describe
the twistor space $J(M,\om)$ for manifolds of
Ricci-type in more detail. In particular, since manifolds of Ricci-type
can always be obtained locally by the reduction process explained in
Theorem~\ref{thm:Ricci model}, it seems reasonable to expect the bundle
$J(M, \om)$ and its complex structure to arise as a reduction of
some sort as well. It is the aim of this section to look at this question.

First of all, consider the twistor space $J(\R^{2n+2}, \Om')$ of
$\R^{2n+2}$ with the canonical symplectic structure $\Om'$, let $A \in
\sp(\R^{2n+2}, \Om')$ and $\Sigma_A:= \{ x \in \R^{2n+2} \mid \Om'(x,
Ax) = 1\}$. We decompose the tangent spaces at $x \in \Sigma_A \subset
\R^{n+2}$ as
\[
T_x \Sigma_A = span(Ax) \oplus {\cal D}_x\ \ \ \mbox{and}\ \ \ T_x
\R^{2n+2} = span(x, Ax) \oplus {\cal D}_x,
\]
where ${\cal D}_x$ is determined by $\Om'(span(x, Ax), {\cal D}_x) = 0$.
Moreover, we have the canonical projection $\pi: \Sigma_A \supset U \ra
M^{red}$ where $U$ is a sufficiently small open subset such that
$M^{red} = U/(\exp \R A)^{loc}$ is a manifold, as explained in
Section~\ref{section:reduction}. Recall that the symplectic structure
$\om$ on $M^{red}$ is determined by the requirement that $d\pi|_{{\cal
D}_x}: ({\cal D}_x, \Om') \ra (T_{\pi(x)} M^{red}, \om_{\pi(x)})$
becomes a symplectic isomorphism. Moreover, the Ricci-type connection
$\nabla$ on $(M^{red}, \om)$ is defined by
\[
\ov{\nabla_X Y} = \left(\nabla_{\ov X}^0 \ov Y\right)_{\cal D},\ \
\mbox{for all vector fields $X, Y$ on $M^{red}$},
\]
where the bar denotes the horizontal lift w.r.t.\ the distribution ${\cal
D}$, the subscript ${\cal D}$ denotes the image under the canonical
projection $T_x \R^{2n+2} = span(x, Ax) \oplus {\cal D}_x \ra {\cal
D}_x$, and where $\nabla^0$ denotes the canonical connection on
$\R^{2n+2}$.

\begin{lemma}
Let $p: J(M^{red}, \om) \ra M^{red}$ be the twistor fibration. Then the
pull-back of $p$ under the map $\pi: \Sigma_A \ra M^{red}$ is given by
the fibration $\hat p: \hat Z_A \ra \Sigma_A$ where
\[
\hat {Z}_A:= \left\{ (x,J) \in J(\R^{2n+2}, \Om') \mid x \in \Sigma_A,
Jx = Ax,\ \ JAx = -x \right\},\quad \hat p(x,J) = x
\]
and the map $\hat \pi : \hat Z_A \ra J(M^{red}, \om)$ is given by
\[
(\hat \pi (J)) (v):= d\pi (J \ov v), \quad\forall v \in T_{\pi(x)}
M^{red}.
\]
Moreover, $J(M^{red}, \om)$ is the quotient of $\hat Z_A$ under the (local) 
action of $\exp (\R A)$ given by
\[
\exp(t A) \cdot (x,J):= (\exp (tA)x, \exp (tA)\,J\,\exp(-tA)).
\]
\end{lemma}

\begin{proof}
Let $x \in \Sigma_A$ and consider an element $J \in j(T_{\pi(x)}
M^{red}, \om)$. Using the symplectic isomorphism $d\pi_x: ({\cal D}_x,
\Om') \ra (T_{\pi(x)} M^{red}, \om)$, we see that there is a unique way
to extend from ${\cal D}_x$  to $\R^{2n+2}$ to give an element $\hat J 
\in j(\R^{2n+2},\Om')$ with $(x,\hat J) \in \hat Z_A$ and $\hat\pi
(x,\hat J) = J$. This shows the first assertion. The second assertion is
straightforward.
\end{proof}

After having described $J(M^{red}, \om)$ as a manifold as the reduction
of $\hat Z_A$ under the (local) action of $\exp(\R A)$ it would have
been nice to see that the twistor almost complex structure on
$J(M^{red}, \om)$ coming from the reduced connection arose as a quotient
from the twistor complex structure on $J(\R^{2n+2}, \Om')$ coming from
$\nabla^0$. Unfortunately this is not the case. The inverse image by
$d\hat\pi$ of the twistor almost complex structure on $J(M^{red}, \om)$
gives an almost complex structure on the quotient of $T\hat Z_A$ by $\R
\hat H_A$ and this would be induced from $J^{\nabla^0}$ on $J(\R^{2n+2},
\Om')$ only if $J^{\nabla^0} T\hat Z_A \subset T\hat Z_A + \R
J^{\nabla^0}\hat H_A$. A straightforward calculation shows that this
last condition is equivalent to
\[
[A,J] {\cal D}_x = 0, \quad \forall (x,J) \in \hat Z_A
\]
and that this condition is too strong for general $A$.

\section{Ricci-flat connections}\label{sec:Ricciflat}

\subsection{A construction by induction}\label{section:induction}

  \begin{definition}
A {\bf{contact quadruple}} $(M, N,\alpha, \pi)$ is a $2n$
dimensional smooth  manifold $M$, a $2n+1$ dimensional smooth manifold
$N$, a co-oriented contact structure $\alpha$ on $N$ (i.e. $\alpha$
is a $1$-form on $N$ such that $\alpha\wedge (d\alpha)^n$ is nowhere
vanishing), and a smooth submersion $\pi: N\rightarrow M$ with $d\alpha
=\pi^*\omega$ where $\omega$ is a symplectic $2$-form on $M$.
 \end{definition}
 \begin{definition}
Given a contact quadruple $(M, N,\alpha, \pi)$ the {\bf{induced
symplectic manifold}} is the $2n+2$ dimensional manifold
$$
P:=N\times \R$$ endowed with the (exact) symplectic structure
$$
\mu:=2e^{2s}~ds\wedge p_1^*\alpha +e^{2s}~ dp_1^*\alpha =d(e^{2s}~p_1^*\alpha)
$$
 where $s$ denotes the variable along $\R$ and  
$p_1: P \rightarrow N$ the projection on the first factor. 
\end{definition}
Induction in the sense of building a $(2n+2)$-dimensional symplectic manifold
from a symplectic manifold of dimension $2n$ is also considered
by Kostant in \cite{Kostant}.

\begin{remark}~

\noindent$\bullet$ The vector field $S:=\partial_s$ on $P$ is such that
$i(S)\mu=2{\rm e}^{2s}(p_1^*\alpha)$; hence 
$L_{S}\mu=2\mu$ and $S$ is a conformal vector
field.

\noindent$\bullet$ The Reeb vector field $Z$ on $N$ (i.e. the vector
field $Z$ on $N$ such that $i(Z)d\alpha=0$ and $i(Z)\alpha=1$) lifts to
a vector field $E$ on $P$ such that: $p_{1*}E=Z$ and $ds(E)=0.$ Since
$i(E)\mu=-d({\rm e}^{2s})$, $E$ is a Hamiltonian vector field on
$(P,\mu)$.  Furthermore
\begin{eqnarray*}
   [E,S] &=& 0 \\
   \mu(E,S) &=& -2{\rm e}^{2s}.
\end{eqnarray*}

\noindent$\bullet$ Observe also that if $\Sigma=\{\, y\in P\, |\,
s(y)=0\, \}$, the reduction of $(P,\mu$) relative to the constraint
manifold $\Sigma$ (which is isomorphic to $N$) is precisely
$(M,\omega)$.
 
\noindent$\bullet$ For $y\in P$ define $H_y(\subset T_yP)=
>E,S<^{\bot_\mu}$. Then $H_y$ is symplectic and $(\pi\circ p_1)_{*y}$
defines a linear isomorphism between $H_y$ and $T_{\pi p_1(y)}M$. Vector
fields on $M$ thus admit ``horizontal'' lifts to $P$.
 \end{remark}
We shall now prove that any symplectic connection $\nabla$ on
$(M,\omega)$ can be lifted  to a symplectic connection
on $(P,\omega^P)$ which is {\bf Ricci-flat}.
We shall initially  define a connection $\nabla^P$ on $P$ induced by $\nabla$.

First some notation:\\
$p$ denotes the projection $p=\pi\circ p_1:P\rightarrow M$.\\
If $X$ is a vector field on $M$, ${\bar{\bar{X}}}$
 is the vector field on $P$ such that
$$
(i)~~p_*{\bar{\bar{X}}}=X \quad\quad  (ii)~~(p_1^*\alpha)({\bar{\bar{X}}})=0
\quad\quad  (iii)~~ds({\bar{\bar{X}}})=0.
$$
Recall that $E$ is the vector field on $P$ such that
$$
(i)~~p_{1*}E=Z \quad\quad  (ii)~~ds(E)=0.
$$
Clearly  the values at any point of $P$ of the vector fields 
 ${\bar{\bar{X}}},E, S=\partial_s$ span the tangent space to 
$P$ at that point and we have 
$$
[E,\partial_s]=0 \quad
[E,{\bar{\bar{X}}}]=0æ\quad
[\partial_s,{\bar{\bar{X}}}]=0 \quad
[{\bar{\bar{X}}},{\bar{\bar{Y}}}]={\overline{\overline{[X,Y]}}}-p^*\omega(X,Y)E.
$$
The formulas for $\nabla^P$ are:
\begin{eqnarray*}
\nabla^P_{{\bar{\bar{X}}}} {\bar{\bar{Y}}}&=&\overline{\overline{{\nabla}_X Y}}
   -\half p^*(\omega(X,Y))E- p^*({\hat{s}}(X,Y))\partial_s
      \\[2mm]
\nabla^P_{E} {\bar{\bar{X}}}&=&\nabla^P_{{\bar{\bar{X}}}}E =2
    \overline{\overline{{\sigma}X}} + p^*(\omega(X,u))\partial_s
      \\[2mm]
\nabla^P_{\partial_s}{\bar{\bar{X}}}&=&
\nabla^P_{{\bar{\bar{X}}}}{\partial_s}=\overline{\overline{X}}
      \\[2mm]
\nabla^P_E E&=&p^*f \,{\partial_s}-2\overline{\overline{U}} 
      \\[2mm]
\nabla^P_E{\partial_s}&=&\nabla^P_{\partial_s}E=E
      \\[2mm]
\nabla^P_{\partial_s}{\partial_s}&=&{\partial_s} 
\end{eqnarray*}
where  $f$ is a function on $M$, $U$ is a vector field on $M$, 
${\hat{s}}$ is a symmetric $2$-tensor on
$M$, and $\sigma$ is the endomorphism of $TM$ associated to $s$, hence
${\hat{s}}(X,Y)=\omega(X,\sigma Y)$.

These formulas have the correct linearity properties
and yield a torsion free linear connection on $P$.
One checks readily that $\nabla^P\mu=0$ so that $\nabla^P$
is a symplectic connection on $(P,\mu)$.

We now compute the curvature $R^{\nabla^P}$ of this connection $\nabla^P$. 
We get
\begin{eqnarray*}
R^{\nabla^P}({\bar{\bar{X}}},{\bar{\bar{Y}}}){\bar{\bar{Z}}}&=&
\overline{\overline{R^{\nabla}(X,Y)Z}}\cr
&&+\overline{\overline{2\omega(X,Y)\sigma Z-\omega(Y,Z)\sigma X
+\omega(X,Z)\sigma Y-{\hat{s}}(Y,Z)X+{\hat{s}}(X,Z)Y}} \cr
&&+p^*[\omega(X,D(\sigma,U)(Y,Z))
-\omega(Y,D(\sigma,U)(X,Z)] \partial_s\cr
&& ~~\cr
R^{\nabla^P}({\bar{\bar{X}}},{\bar{\bar{Y}}})E&=&
\overline{\overline{2D(\sigma,U)(X,Y)-2D(\sigma,U)(Y,X)}}\cr
&&+p^*[\omega(X,\half fY-\nabla_YU-2{\sigma}^2Y)
-\omega(Y,\half fX-\nabla_XU-2{\sigma}^2X)]\partial_s\cr
&& ~~\cr
R^{\nabla^P}({\bar{\bar{X}}},E){\bar{\bar{Y}}}&=&
\overline{\overline{2D(\sigma,U)(X,Y)}}-
p^*[\omega(Y,\half fX-\nabla_XU-2{\sigma}^2X) ] \partial_s\cr
&& ~~\cr
R^{\nabla^P}({\bar{\bar{X}}},E)E&=&
2\overline{\overline{\half fX-\nabla_XU-2{\sigma}^2X}}
+p^*[Xf+4s(X,u)] \partial_s\cr
&& ~~\cr
R^{\nabla^P}({\bar{\bar{X}}},{\bar{\bar{Y}}}){\partial_{s}}&=&0 \quad \quad
R^{\nabla^P}({\bar{\bar{X}}},E){\partial_{s}}=0 \cr
R^{\nabla^P}({\bar{\bar{X}}},{\partial_{s}}){\bar{\bar{Y}}}&=&0 \quad\quad
R^{\nabla^P}({\bar{\bar{X}}},{\partial_{s}})E=0 \quad 
R^{\nabla^P}({\bar{\bar{X}}},{\partial_{s}}){\partial_{s}}=0\cr
R^{\nabla^P}(E,{\partial_{s}}){\bar{\bar{X}}}&=&0\quad \quad
R^{\nabla^P}(E,{\partial_{s}})E=0\quad
R^{\nabla^P}(E,{\partial_{s}}){\partial_{s}}=0
\end{eqnarray*}
where
$$
D(\sigma,U)(Y,Y'):=(\nabla_Y\sigma)Y'+\half \omega(Y',U)Y-\half\omega(Y,Y')U.
$$
The Ricci tensor $r^{\nabla^P}$ of the connection $\nabla^P$ is given by
\begin{eqnarray*}
r^{\nabla^P}({\bar{\bar{X}}},{\bar{\bar{Y}}})&=&
              r^{\nabla}(X,Y)+ 2(n+1) {\hat{s}}(X,Y)\cr
r^{\nabla^P}({\bar{\bar{X}}},E)&=&
            -(2n+1)\omega(X,u) -2\Tr[Y\rightarrow (\nabla_Y\sigma)(X)]\cr
r^{\nabla^P}({\bar{\bar{X}}},{\partial_s})&=&0\cr
r^{\nabla^P}(E,E)&=&4\Tr (\sigma^2)-2nf
            +2\Tr[X\rightarrow \nabla_XU ] \cr
r^{\nabla^P}(E,{\partial_{s}})&=&0\cr
r^{\nabla^P}({\partial_{s}},{\partial_{s}})&=&0
\end{eqnarray*}

\begin{theorem}\cite{CahenGutt}\label{th:Ricciflat}
In the framework described above, $\nabla^P$ is a symplectic connection
on $(P,\mu)$ for any choice of ${\hat{s}},U$ and $f$.
The vector field $E$ on $P$ is affine 
( $L_{\tilde{E}}\nabla^P=0$) and symplectic ( $L_{\tilde{E}}\mu=0$);
the vector field $\partial_s$ on $P$ is  affine and conformal 
($L_{\partial_s}\mu=2\mu$).

Furthermore, choosing
\begin{eqnarray*}
{\hat{s}}&=&\frac{-1}{2(n+1)}r^\nabla\cr
{\underline{U}}:
&=&\omega(U,\cdot)=\frac{2}{2n+1}\Tr[Y\rightarrow\nabla_Y\sigma]\cr
f&=&\frac{1}{2n(n+1)^2}\Tr (\rho^\nabla)^2 +\frac{1}{n}
\Tr [X\rightarrow \nabla_XU].
\end{eqnarray*}
we have:
\begin{itemize}
\item  the connection $\nabla^P$ on 
$(P,\mu)$ is Ricci-flat (i.e. has zero 
Ricci tensor); 
\item if the symplectic connection $\nabla$ on $(M,\omega)$ is of
Ricci-type, then 
the connection $\nabla^P$ on $(P,\mu)$ is flat.
\item if the connection $\nabla^P$ is locally symmetric, the connection
$\nabla$ is of Ricci-type, hence  $\nabla^P$ is flat.
\end{itemize}
\end{theorem}
\begin{proof}
The first point is an immediate consequences of the formulas above for
$r^{\nabla^P}$. The second point is a consequence of the differential
identities satisfied by the Ricci-type symplectic connections. The third
point comes from the fact that $(\nabla^P_{\bar{\bar{Z}}} R^{\nabla^P})
({\bar{\bar{X}}},{\bar{\bar{Y}}}){\bar{\bar{T}}}$ contains only one term
in $E$ whose coefficient is $\half W^{\nabla^P}(X,Y,T,Z)$.
\end{proof}

\subsection{Examples of contact quadruples}
We give here examples of contact quadruples
corresponding to a given symplectic manifold $(M,\omega)$
(i.e. examples of $(N,\alpha, \pi)$
 where $N$ is a smooth $2n+1$ dimensional manifold, $\alpha$ 
 is a $1$-form on $N$ such that
 $\alpha\wedge (d\alpha)^n$ is nowhere vanishing, $\pi: N\rightarrow M$
 is a smooth submersion and $d\alpha =\pi^*\omega$.
\begin{itemize}
\item  Let $(M, \omega=d\lambda)$ be an {\bf exact symplectic  manifold}.
Define $N=M\times\mathbb R$, $\pi=p_1$
(=projection of the  first factor), $\alpha=dt+p_1^*\lambda$;
then $(N,\alpha)$ is a contact  manifold and $(M,N,\alpha,\pi)$
is a contact quadruple.

The associated induced manifold is
$P=N\times\mathbb R=M\times\mathbb R^2$; with coordinates 
$(t,s)$ on $\R^2$ and obvious identification
$$  \mu={\rm e}^{2s}~[\,d\lambda +
2ds\wedge  (dt+\lambda)\,].
$$

\item Let $(M,\omega)$ be a {\bf quantizable symplectic manifold}; this
means that there is a complex line bundle $L\map{p}M$ with hermitean
structure $h$ and a connection $\nabla$ on $L$ preserving $h$ whose
curvature is proportional to $i\omega$.\\
Define $N: = \{\,\xi\in L~|~h(\xi,\xi\,)=1\}\subset L$ to be the
unit circle sub-bundle. It is a principal  $U(1)$ bundle and $L$ is
the associated bundle $L=N\times_{U(1)}\mathbb C$. The connection
1-form on $N$ (representing $\nabla$) is $u(1)=i\R$ valued
and will be denoted $\alpha'$; its curvature
is $d\alpha'=ik\omega$. Define $\alpha:=\dfrac{1}{ik}\alpha'$ 
and $\pi:=p|_N: N\rightarrow M$  the surjective submersion.
Then $(M,N,\alpha,\pi)$ is a contact quadruple.

The associated induced manifold $P$ is in bijection with
$L_0=L\setminus$ zero section.

\item 
 Let $(M,\omega)$ be a {\bf connected homogeneous symplectic
 manifold}; i. e. $M=G/H$ where $G$ is a Lie group which we may
assume  connected and simply connected and where $H$ is the
stabilizer in  $G$ of a point $x_0\in M$. If $p: G\rightarrow
M:g\rightarrow gx_0,\quad p^*\omega$ is a left invariant closed 2-form
on $G$ and  $\Omega=(p^*\omega)_e$, (e=neutral element of $G$) is a
Chevalley  2-cocycle on ${\g}$ (=Lie Algebra of $G$) with values
in $\mathbb R$ (for the trivial representation).\\
Notice that $\Omega$ vanishes as soon as one  of its arguments is
 in ${\h}$ (=Lie algebra of $H$). Let ${\g}_1=\g
\oplus\mathbb R$ be the central extension of ${\g}$ defined by
 $\Omega$; i. e.
 $$
  [(X,a),(Y,b)]=([X,Y],\Omega (X,Y)).
  $$
 Let ${\h}'$ be the subalgebra of $\g_1$, isomorphic to $\h$,  
defined by ${\h}':=\{\,(X,0)\,\vert\,X\in \h\,\}$. Let $G_1$
be the connected and simply connected group of algebra $\g$, 
and let $H'$ be the connected subgroup of $G_1$
with Lie algebra ${\h}'$. {\bf{Assume $H'$ is closed}}.\\
Then $G_1/H'$ admits a natural structure of smooth manifold;
define $N:=G_1/H'$.
Let $p_1:G_1\rightarrow G$
be the homomorphism whose differential is the  projection
$\g_1\rightarrow \g$ on the first factor; clearly $p_1(H')\subset H$.
Define $\pi:N=G_1/H'\rightarrow M=G/H: g_1H'\mapsto p_1(g_1)H$; it is a
surjective submersion.

The contact form $\alpha$ on $N$ is constructed as follows: $p_1^*\circ
p^*\omega$ is a left invariant  closed 2-form on $G_1$ vanishing on the
fibres of $p\circ p_1:G_1\rightarrow M$. Its  value $\Omega_1$ at the
neutral element $e_1$ of $G_1$ is a  Chevalley 2-cocycle of $\g_1$ with
values in $\mathbb R$. Define the $1$-cochain  $\alpha_1:
\g_1\rightarrow \mathbb R: (X,a)\rightarrow -a$. Then
$\Omega_1=\delta\alpha_1$ is a coboundary. Let $\tilde\alpha_1$ be the
left invariant 1-form on $G_1$ corresponding to $\alpha_1$. Let $q:
G_1\rightarrow G_1/H'=N$ be the natural projection.\\
There exists a $1$-form $\alpha$ on $N$
 so that $q^*\alpha=\tilde\alpha_1$.
Indeed, for any $X\in \h'$ we have
\begin{eqnarray*}
   &&i(\tilde X)\tilde\alpha_1 = \alpha_1(X)=0 \\
   &&(L_{\tilde X}\tilde\alpha_1)(\widetilde (Y,b)) =
-\tilde\alpha_1([\tilde X,\widetilde (Y,b)])=-\alpha_1([X,(Y,b)])=\Omega(X,Y)=0 
\end{eqnarray*}
where $\tilde U$  is the  left invariant vector field on $G_1$
corresponding to $U\in \g_1$. Furthermore $d\alpha=\pi^*\omega$ because
both are $G_1$ invariant $2$-forms on $N$ and:
\begin{eqnarray*}
 (d\alpha)_{q(e_1)}((X,a)^{*N},(Y,b)^{*N})&=
 &(q^*d\alpha)_{e_1}(\widetilde{(X,a)},\widetilde{(Y,b)})\\
 &=&(d\tilde\alpha_1)_{e_1}(\widetilde{(X,a)},\widetilde{(Y,b)})\\
 &=&\Omega(X,Y),\\
 &=&\omega_{x_0}(X^{*M},Y^{*M})\\
 &=&(\pi^*\omega)_{q(e_1)}((X,a)^{*N},(Y,b)^{*N})
  \end{eqnarray*}
where we denote  by $U^{*N}$  the fundamental vector field on $N$
associated to $U \in \g_1$  .

\item
If $(M,\omega,\nabla)$ is a {\bf simply connected symplectic manifold
with a Ricci-type connection} we have seen in Section \ref{globalRicci} 
how to build the manifold $N$ as a holonomy bundle over $M$ corresponding 
to a connection built on the extension $B'(M)$
of the frame bundle $B(M)$ over $M$.
\end{itemize}

\subsection{More about reduction}

Let $(P,\omega^P)$ be a symplectic manifold of dimension $(2n+2)$.
Assume $P$ admits a   conformal vector field $S$:
$$
L_S\omega^P=2\omega^P;\qquad{\rm{define~~}}\alpha:=\half
i(S)\omega^P\qquad {\rm{so~that~}} d\alpha=\omega^P.
$$
Assume also that $P$ admits a symplectic vector field ${\tilde{E}}$
commuting with $S$
$$
L_{\tilde{E}}\omega^P=0\qquad\qquad [S,{\tilde{E}}]=0 
\qquad\qquad (\Rightarrow L_{\tilde{E}}\alpha=0).
$$
Define
$$
\Sigma=\{x\in P\,|\,\omega^P_x(S,{\tilde{E}})=1\}
$$
and assume that it is non-empty. 
The tangent space to the hypersurface $\Sigma$ is given
by 
$$
T_x\Sigma= \ker (i({\tilde{E}})\omega^P)_x={\tilde{E}}^{\perp_{\omega^P}}.
$$
The restriction of
$\omega^P_x$ to $T_x\Sigma$ has rank $2n-2$ and a radical spanned by
${\tilde{E}}_x$.\\
The restriction of $\alpha$ to $\Sigma$ is a contact
$1$-form on $\Sigma$.

Let $\sim$ be the equivalence relation defined on $\Sigma$ by the
flow of ${\tilde{E}}$. Assume that the quotient $\Sigma/\sim$
has a $2n$ dimensional manifold $M$ structure 
so that $\pi:\Sigma\rightarrow\Sigma/\sim=M$
is a smooth submersion.\\
Define on $\Sigma$ a ``horizontal'' distribution of dimension $2n$,
$\mathcal H$, by
$$
\mathcal H=>{\tilde{E}},S<^{\perp_{\omega^P}},
$$
and remark that
$\pi_{*\vert_{{\mathcal{H}}_y}}:{{\mathcal{H}}_y}\rightarrow
T_{x=\pi(y)}M$ is an isomorphism.

Define as usual the reduced 2-form $\omega^M$ on $M$ by
$$
\omega^M_{x=\pi(y)}(Y_1,Y_2)=\omega^P_y(\bar Y_1,\bar Y_2)
$$
where $\bar Y_i$ $(i=1,2)$ is defined by (i) $\pi_*\bar Y_i=Y_i$
(ii) $\bar Y_i\in \H_y$.\\
Notice that $\pi_*[{\tilde{E}},\bar Y]=0$,
and $\omega^P(S,[{\tilde{E}},\bar Y])=-L_{\tilde{E}}\omega^P(S,\bar Y)
+{\tilde{E}}\omega^P(S,\bar Y)=0$  hence 
$$[{\tilde{E}},\bar Y]=0.$$
The definition of $\omega^M_x$ does not depend on the choice of $y$.
Indeed
$$
{\tilde{E}}\omega^P(\bar Y_1,\bar Y_2)=L_{\tilde{E}}\omega^P(\bar Y_1,\bar Y_2)
+\omega^P([{\tilde{E}},\bar
Y_1],\bar Y_2)+\omega^P(\bar Y_1,[{\tilde{E}},\bar Y_2])=0.
$$
Clearly $\omega^M$ is of maximal rank $2n$ as $\mathcal H$ is a symplectic
subspace. Finally
\begin{eqnarray*}
  \pi^*(d\omega^M(Y_1,Y_2,Y_3)) &=& 
  \cyclic_{123} ~(Y_1\omega^M(Y_2,Y_3)-\omega^M([Y_1,Y_2],Y_3) )\\
   &=& \cyclic_{123} ~ (\bar Y_1\omega^P(\bar Y_2,\bar
   Y_3)-\omega^P(\overline{[Y_1,Y_2]}, \overline{Y}_3))
\end{eqnarray*}
and
$$
[\bar Y_1,\bar Y_2]=\overline{[Y_1,Y_2]}+\omega^P(S,[\bar Y_1,\bar
Y_2]){\tilde{E}}.
$$
Hence $\omega^M$ is closed and thus symplectic. Clearly 
$\pi^*\omega^M=\omega^P_{\vertÑ\Sigma}=d(\alpha_{\vertÑ\Sigma})$.

\begin{remark}  The manifold $M$ is the first element of
a contact quadruple $(M,\Sigma,\half\alpha_{\vert_\Sigma},\pi)$.
\end{remark}

We shall now consider the reduction of a connection.
Let $(P,\omega^P),{\tilde{E}},S,\Sigma, M,\omega^M$ be as above.
Let $\nabla^P$ be a symplectic connection on $P$ and assume that
the vector field ${\tilde{E}}$ is affine ($L_{\tilde{E}}\nabla^P=0$).

Then define a connection $\nabla^{\Sigma}$ on $\Sigma$ by
$$
\nabla_A^{\Sigma}B:=\nabla_A^PB-\omega^P(\nabla_A^PB,{\tilde{E}})S
=\nabla_A^PB+\omega^P(B,\nabla_A^P{\tilde{E}})S.
$$
Then $\nabla^\Sigma$ is a torsion free connection and ${\tilde{E}}$ is an
affine vector field for $\nabla^\Sigma$.

 Define a connection
$\nabla^M$ on $M$ by:
$$
\overline{\nabla_{Y_1}^MY_2}(y)=\nabla_{\bar Y_1}^\Sigma\bar
Y_2(y)-\omega^P(\bar Y_2,\nabla_{\bar Y_1}^P S){\tilde{E}}.
$$
If $x\in M$, this definition does not depend on the choice of $y\in
\pi^{-1}(x)$ and one can check that the connection $\nabla^M$ is symplectic.
\begin{lemma}\cite{CahenGutt}
Let $(P,\omega^P)$ be a symplectic manifold admitting a symplectic
connection $\nabla^P$, a conformal vector field $S$, a symplectic vector
field ${\tilde{E}}$ which is affine and commutes with $S$. If the
constraint manifold $\Sigma=\{x\in P|\omega_x(S,{\tilde{E}})=1\}$ is not
empty, and if the reduction of $\Sigma$ is a manifold $M$, this manifold
admits a symplectic structure $\omega^M$ and a natural reduced
symplectic connection $\nabla^M$.
\end{lemma}
In particular
\begin{theorem}\cite{CahenGutt}
Let $(P,\omega^P)$ be a symplectic manifold admitting a conformal vector
field $S$ ($L_S\mu=2\mu$) which is complete,  a symplectic vector field
${\tilde{E}}$ which  commutes with $S$ and assume that, for any $x\in
P,~\mu_x(S,{\tilde{E}})>0$. If the reduction of $\Sigma=\{x\in
P\,|\,\mu_x(S,{\tilde{E}})=1\}$ by the flow of ${\tilde{E}}$ has a
manifold structure $M$ with $\pi:\Sigma\rightarrow M$ a surjective
submersion, then $M$ admits a reduced symplectic structure $\omega^M$
and $(P,\omega^P)$ is obtained by induction from $(M,\omega^M)$ using
the contact quadruple $(M,\Sigma, \half
i(S)\omega^P_{\vert_\Sigma},\pi)$.

In particular $(P,\omega^P)$ admits a Ricci-flat connection.
 \end{theorem}

\section{Non-commutative symplectic symmetric spaces}
\subsection{Motivations}

After the celebrated example of the quantum torus and related
non-commutative spaces \cite{Rieffel,CoLa}, it appeared natural to try
to define non-commutative spaces through oscillatory integral formulae
in the same spirit of \cite{Rieffel} but with larger symmetry
groups---other  than $\R^d$--- implementing this way in our class of
non-commutative manifolds not only an operator algebraic framework but
also a strong geometric content.

\noindent In the context of symmetric spaces, this leads to the
so-called `WKB-quantisation of symplectic symmetric spaces' as
introduced by Karas\"ev, Weinstein and Zakrzewski in the mid 90's (see
\cite{Wei} and references therein). Originally (cf. \cite{Wei}), the
following  question was raised in the framework of Hermitean symmetric
spaces of the non-compact type $M=G/K$.
\begin{question}
Given a $G$-invariant product on $M$ expressed in
 the `WKB form':
\begin{equation}\label{PROD}
(u\star_\theta v)(x)=\frac{1}{\theta^{2n}}\int_{M\times
M}a_\theta(x,y,z)\, \exp\left( \frac{i}{\theta}\,S(x,y,z)\right)\,
u(y)\, v(z) \, dy\, dz\;,
\end{equation}
which conditions on the phase function $S\in C^\infty(M^3,\R) $ and on
the amplitude $a_\theta=a_0+\theta a_1+\theta^2 a_2+...\in
C^\infty(M^3,\R)[[\theta]]$ does one need in order to guarantee formal
associativity of $\star_\theta$?
\end{question}

Weinstein provided some evidence showing that the phase function $S$
should be closely related to the three point function denoted hereafter
$S_{\mathrm{can}}$ and (partially) defined as follows. Given three
points $x,y$ and $z$ in the symmetric space with the property that the
equation $s_xs_ys_z(X)=X$ admits a (unique) solution $X$, the value of
$S_{\mathrm{can}}(x,y,z)$ is given by the symplectic area
$$
S_{\mathrm{can}}(x,y,z)=\int_{\stackrel{\Delta}{XYZ}}\omega,
$$
where $\stackrel{\Delta}{XYZ}$ denotes the geodesic triangle with vertices
$$
X,\qquad Y:=s_z(X),\qquad Z=s_ys_z(X).
$$
The question of characterising geometrically the amplitude was
left open.

Another important aspect in this problematic is to determine whether
such a WKB product underlies topological function algebras analogous to
the continuous field of $C^\star$-algebras deforming $C(\TT)$ in the
case of the quantum torus $\TT_\theta$. More precisely:

\begin{question}\label{FA}Given 
$\theta$ in some deformation parameter space, does one have a function
space $\A_\theta$
$$
C^\infty_c(M)\subset \A_\theta\subset{\cal D}'(M)
$$
such that the pair $(\A_\theta,\star_\theta)$ is a topological 
$G$-algebra\footnote{that is, an associative algebra underlying a 
topological vector space on which the group $G$ acts strongly continuously
by algebra automorphisms.}?
\end{question}
In the present section, we survey a geometrical approach to these
questions initiated in \cite{Biel2} and based on the definition of a
class of three point functions, called hereafter `admissible', on
symplectic symmetric spaces. Admissible functions are characterised by
compatibility properties with the symmetries of the symmetric space at
hand. We will show how these properties guarantee associativity of the
oscillating product associated to such an admissible function in the
case of a cocyclic function. We will then indicate on a curved two
dimensional example (a generic coadjoint orbit of the Poincar\'e group
in dimension 1+1)  how the cocycle condition can be relaxed by
introducing a non-trivial amplitude in the oscillating kernel. We will
end by mentioning a result which solves Question \ref{FA} above in this
particular context. Let us add that the general results obtained in this
direction go far beyond the particular example presented here. For
instance, the solution for the general solvable symmetric case has led
to several universal deformation formulae for actions of various classes
of solvable Lie groups \cite{BMas,BMae,BBM}. Applications in
non-commutative geometry (non-commutative causal black holes
\cite{BDRS}) as well as in analytic number theory (Rankin--Cohen brackets
on modular forms \cite{BXY}) have been developed.

\subsection{Basic definitions and the cocyclic case}
Denoting by $K:=a_\theta\;e^{\frac{i}{\theta}S}$ the oscillating kernel
defining the product given in (\ref{PROD}), a simple computation shows
that associativity of $\star_\theta$ is (at least formally) equivalent
to the following condition:

\begin{equation}\label{ASSK}
\int_{M}K(a,b,t)K(t,c,d)\mu(t)=\int_{M}K(a,\tau,d)K(\tau,b,c)\mu(\tau),
\end{equation}
for every quadruple of points $a,b,c,d$ in $M$. Equation (\ref{ASSK})
obviously holds if one can pass from one integrand to the other using a
change of variables $\tau=\varphi(t)$. This motivates

\begin{definition}\label{coucou}
Let $(M,\mu)$ be an orientable manifold endowed with a volume form
$\mu$. A three-point kernel $K\in C^{\infty}(M\times M\times M)$ is {\bf
geometrically associative} if for every quadruple of points $a,b,c,d$ in
$M$ there exists a volume preserving diffeomorphism
$$
\varphi:(M,\mu)\to(M,\mu),
$$
such that for all $t$ in $M$:
$$
K(a,b,t)K(t,c,d)=K(a,\varphi(t),d)K(\varphi(t),b,c).
$$
\end{definition}
In the sequel, we give sufficient conditions  for geometric
associativity. From now on, $(M,\omega,s)$ denotes a symplectic
symmetric space (not necessarily Hermitean). We first observe the
following group-like cohomological complex naturally associated to our
context.
\begin{definition}
A $k$-cochain on $M$ is a totally skew symmetric real-valued smooth
function $S$ on $M^k$ which is invariant under the (diagonal) action of
the symmetries $\{s_x\}_{x\in M}$. Denoting by $P^k(M)$, the space of
$k$-cochains, one has the cohomology operator $\delta:P^k(M)\to
P^{k+1}(M)$ defined as
\begin{equation}
(\delta S)(x_0,...,x_k):=\sum_j(-1)^jS(x_0,...,\hat{x_j},...,x_k).
\end{equation}
\end{definition}
\begin{definition}\label{ADM}
A  3-cochain $S\in P^3(M)$ is called {\bf admissible} if for all $x \in
M$, one has
$$
S(x,y,z)=-S(x,s_x(y),z)\qquad \forall y,z\in M.
$$
A {\bf Weyl triple} is the data of a symmetric space $M$ endowed with an
invariant volume form $\mu$ together with the data of an admissible
3-cocycle $S$ (i.e. $\delta S=0$).
\end{definition}
Skew-symmetry  naturally leads us to adopt the following 
``oriented graph" type notation for a 3-cochain $S$:

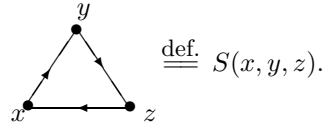
$$\begin{picture}(200,50)
\put(0,0){$x$}
\put(50,0){$z$}
\put(25,40){$y$}
\put(5,5){\trianglehorizontal}
\put(57,19){$\stackrel{\mbox{\small def.}}{=\!=}\,S(x,y,z).$}
\end{picture}$$
A change of orientation in such an ``oriented triangle'' leads to a
change of sign of its value. However, the value represented by such a
``triangle'' does not depend on the way it ``stands'', only the data of
the vertices and the orientation of the edges matters.

Now, consider a Weyl triple $(M,\mu,S)$, and let $A$ be some
(topological) associative algebra. And, for compactly supported
functions $u$ and $v\in C^{\infty}_{c}(M,A)$, consider the following
``product'':
$$
u\star v(x)=\int_{M\times M} u(y) v(z) e^{iS(x,y,z)}\mu(y) \,\mu(z).
$$
\noindent With the above notation for $S$, associativity for $\star$ now
formally takes the following form:

\begin{equation}\label{ASS}
\begin{array}{ccccc}
\begin{picture}(50,100)
\put(0,50){${\displaystyle \int}_M \exp\;i($}
\end{picture} & \carrehoriz{$t$} & 
\begin{picture}(80,100)
\put(0,50){$)\;\mu(t)={\displaystyle \int}_M \exp\;i($}
\end{picture} & \carrevertic{$\tau$} & 
\begin{picture}(50,100)
\put(0,50){$)\;\mu(\tau)\;,$}
\end{picture}
\end{array}
\end{equation}
for every 
quadruple of points $a,b,c,d$ in $M$.

In the above formula, the diagram in the argument of the exponential in
the LHS (respectively the RHS) stands for $S(a,b,t)+S(t,c,d)\;$
(respectively $S(a,d,\tau)+S(\tau,b,c)$).

\begin{prop}\label{GEOM}
Let $(M,\mu,S)$ be a Weyl triple. Then, the associated three-point
kernel $K=e^{iS}$ is geometrically associative.
\end{prop}
\Pf
Fix four points $a,b,c,d$. Regarding \dref{coucou} and 
formula~(\ref{ASS}), one 
needs to construct our volume preserving diffeomorphism 
$\varphi:(M,\mu)\to(M,\mu)$ in such a way that for all $t$,
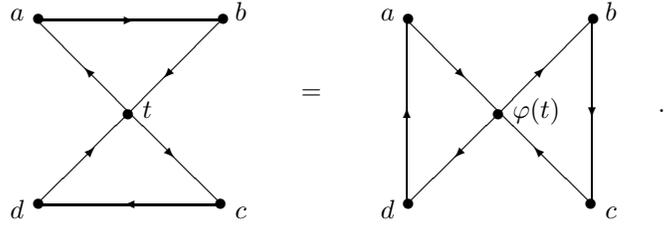
$$
\begin{array}{ccc}
\carrehoriz{$t$} & 
\begin{picture}(20,100)
\put(0,50){=}
\end{picture} & \carrevertic{$\varphi(t)$} 
\end{array}.
$$
We first observe 
that the data of four points $a,b,c,d$ determines what we call
an ``$S$-barycentre", that is a point $g=g(a,b,c,d)$ such that

$$
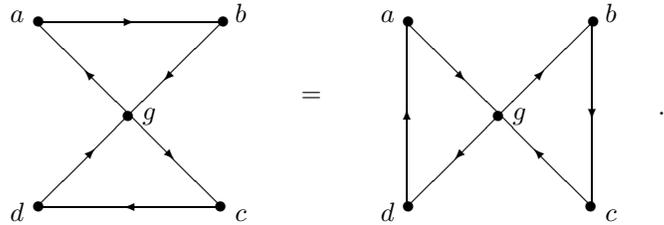
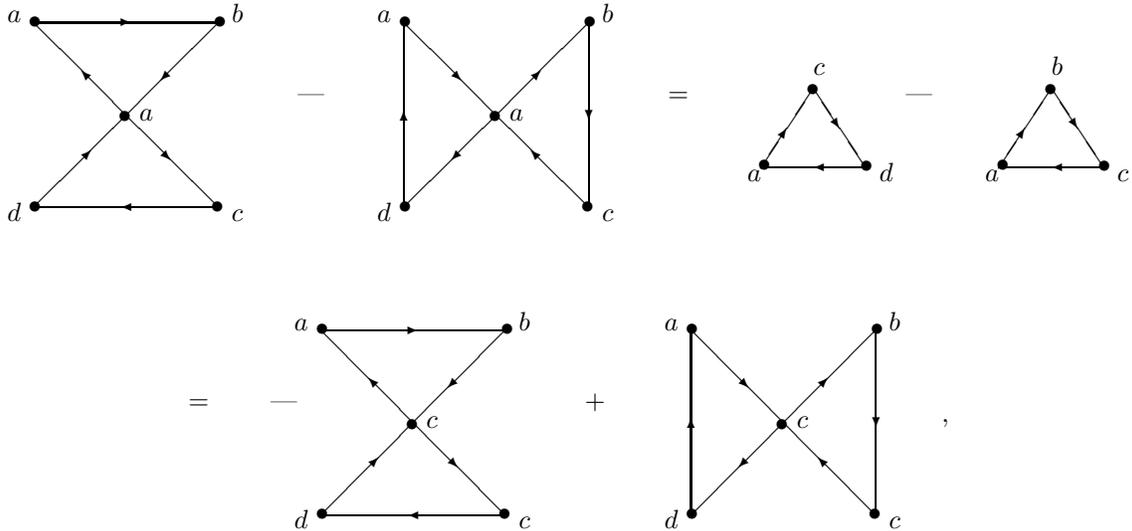
\begin{array}{ccc}
\carrehoriz{$g$} & 
\begin{picture}(20,100)
\put(0,50){=}
\end{picture} & \carrevertic{$g$} 
\end{array}.
$$
Indeed, since
$$
\begin{array}{ccccccc}
\carrehoriz{$a$} &
\begin{picture}(20,100)
\put(0,50){---}
\end{picture} & \carrevertic{$a$} &
\begin{picture}(10,100)
\put(0,50){=}
\end{picture} & \Triangle{a}{c}{d} & \begin{picture}(10,100)
\put(0,50){---}
\end{picture} & \Triangle{a}{b}{c}
\end{array}
$$

$$
\begin{array}{ccccc}
\begin{picture}(30,100)
\put(0,50){= \qquad ---}
\end{picture} &
\carrehoriz{$c$} &
\begin{picture}(20,100)
\put(0,50){+}
\end{picture} & \carrevertic{$c$}
\end{array},
$$
any continuous path joining $a$ to $c$ contains such a point $g$.

Now, we fix once for all such an $S$-barycentre $g$ for $\{ a,b,c,d\}$ and 
we adopt the following notation. For all $x$ and $y$ in $M$, the value 
of $S(g,x,y)$ is denoted by a ``thickened arrow'':
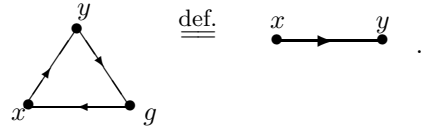
$$
\begin{array}{ccc}
\Triangle{x}{y}{g} & 
\begin{picture}(30,50)
\put(0,50){\mbox{   $\stackrel{\mbox{\small def.}}{=\!=}$}}
\end{picture} & 
\begin{picture}(38,55)
\thicklines
\put(-2,55){$x$}
\put(38,55){$y$}
\put(-2,48){$\bullet$}
\put(38,48){$\bullet$}
\put(20,50){\vector(1,0){2}}
\put(0,50){\line(1,0){40}}
\end{picture}
\end{array}\quad.
$$
Again, a change of orientation in such an arrow changes the sign of its 
value. Also, admissibility  which 
has the form  
$$
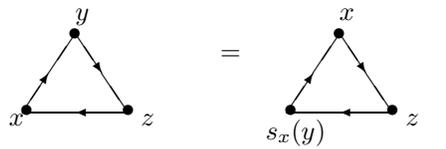
\begin{array}{ccc}
\begin{picture}(60,50)
\put(0,0){$x$}
\put(50,0){$z$}
\put(25,40){$y$}
\put(5,5){\trianglehorizontal}
\end{picture} & 
\begin{picture}(20,50)
\put(10,25){=}
\end{picture} & 
\begin{picture}(60,50)
\put(-3,-5){$s_x(y)$}
\put(50,0){$z$}
\put(25,40){$x$}
\put(5,5){\trianglehorizontal}
\end{picture}
\end{array},
$$
implies
$$
\begin{array}{ccc}
\begin{picture}(50,100)
\thicklines
\put(-2,55){$x$}
\put(38,55){$y$}
\put(-2,48){$\bullet$}
\put(38,48){$\bullet$}
\put(20,50){\vector(1,0){2}}
\put(0,50){\line(1,0){40}}
\end{picture}
&
\begin{picture}(10,90)
\put(0,50){=}
\end{picture} 
&
\begin{picture}(50,100)
\thicklines
\put(-2,55){$x$}
\put(38,55){$s_g(y)$}
\put(-2,48){$\bullet$}
\put(38,48){$\bullet$}
\put(20,50){\vector(-1,0){2}}
\put(0,50){\line(1,0){40}}
\end{picture}
\end{array}
$$
for all $x$ and $y$ in $M$.
While, from cocyclicity, one gets
$$
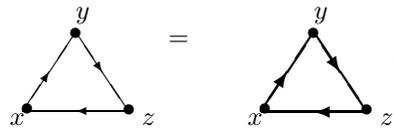
\begin{array}{ccc}
\Triangle{x}{y}{z} &
\begin{picture}(10,90)
\put(0,50){=}
\end{picture} &
\Trianglefat{x}{y}{z}
\end{array}.
$$
Moreover, the barycentric property of $g$ can be written
$$
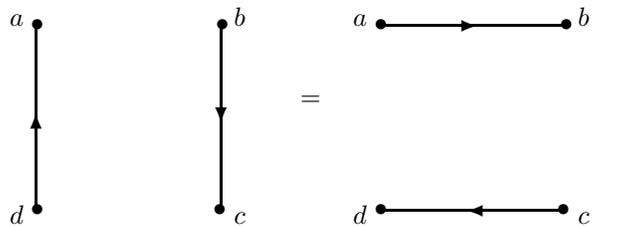
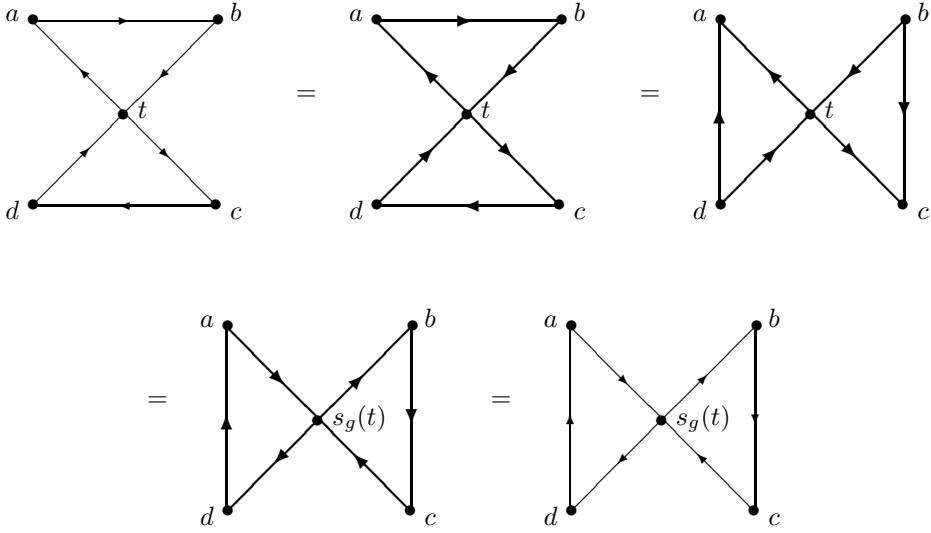
\begin{array}{ccc}
\begin{picture}(100,100)
\thicklines
\put(10,10){\line(0,1){70}}
\put(80,10){\line(0,1){70}}
\put(10,45){\vector(0,1){2}}
\put(80,45){\vector(0,-1){2}}
\put(8,8){$\bullet$}
\put(8,78){$\bullet$}
\put(77,8){$\bullet$}
\put(78,78){$\bullet$}
\put(0,5){$d$}
\put(85,5){$c$}
\put(85,80){$b$}
\put(0,80){$a$}
\end{picture}
&
\begin{picture}(10,90)
\put(0,50){=}
\end{picture} 
&
\begin{picture}(100,100)
\thicklines
\put(10,10){\line(1,0){70}}
\put(10,80){\line(1,0){70}}
\put(45,10){\vector(-1,0){2}}
\put(45,80){\vector(1,0){2}}
\put(8,8){$\bullet$}
\put(8,78){$\bullet$}
\put(77,8){$\bullet$}
\put(78,78){$\bullet$}
\put(0,5){$d$}
\put(85,5){$c$}
\put(85,80){$b$}
\put(0,80){$a$}
\end{picture}.
\end{array}
$$
Hence
$$
\begin{array}{ccccc}
\carrehoriz{$t$}
&
\begin{picture}(10,90)
\put(0,50){=}
\end{picture} 
&
\begin{picture}(100,100)
\thicklines
\put(10,10){\line(1,0){70}}
\put(10,80){\line(1,0){70}}
\put(10,10){\line(1,1){70}}
\put(10,80){\line(1,-1){70}}
\put(30,30){\vector(1,1){2}}
\put(60,60){\vector(-1,-1){2}}
\put(45,10){\vector(-1,0){2}}
\put(45,80){\vector(1,0){2}}
\put(60,30){\vector(1,-1){2}}
\put(30,60){\vector(-1,1){2}}
\put(8,8){$\bullet$}
\put(8,78){$\bullet$}
\put(77,8){$\bullet$}
\put(78,78){$\bullet$}
\put(42,42){$\bullet$}
\put(0,5){$d$}
\put(85,5){$c$}
\put(85,80){$b$}
\put(0,80){$a$}
\put(50,43){$t$}
\end{picture}
&
\begin{picture}(10,90)
\put(0,50){=}
\end{picture} 
&
\begin{picture}(100,100)
\thicklines
\put(10,10){\line(0,1){70}}
\put(80,10){\line(0,1){70}}
\put(10,10){\line(1,1){70}}
\put(10,80){\line(1,-1){70}}
\put(30,30){\vector(1,1){2}}
\put(60,60){\vector(-1,-1){2}}
\put(10,45){\vector(0,1){2}}
\put(80,45){\vector(0,-1){2}}
\put(60,30){\vector(1,-1){2}}
\put(30,60){\vector(-1,1){2}}
\put(8,8){$\bullet$}
\put(8,78){$\bullet$}
\put(77,8){$\bullet$}
\put(78,78){$\bullet$}
\put(42,42){$\bullet$}
\put(0,5){$d$}
\put(85,5){$c$}
\put(85,80){$b$}
\put(0,80){$a$}
\put(50,43){$t$}
\end{picture}
\end{array}
$$

$$
\begin{array}{cccc}
\begin{picture}(10,90)
\put(0,50){=}
\end{picture} 
&
\carreverticfat{$s_g(t)$}
&
\begin{picture}(10,90)
\put(0,50){=}
\end{picture} 
&
\carrevertic{$s_g(t)$}.
\end{array}
$$
One can therefore choose our diffeomorphism $\varphi$ as 
$$
\varphi=s_g.
$$
\EPf

\begin{remark}
Given a symplectic symmetric space Weinstein's function $S_{\mbox{\rm
can}}$ turns out to be admissible (wherever it's well-defined)
\cite{Biel2}. However, the curvature is the obstruction to the
cocyclicity of $S_{\mathrm{can}}$.
\end{remark}

\subsection{A curved example: $SO(1,1)\times\R^2/\R$}
In what follows, we analyse in some details the case of the solvable
symplectic symmetric surface $M=SO(1,1)\times\R^2/\R$. As a homogeneous
symplectic manifold it can be realised as a generic coadjoint orbit of
the Poincar\'e group $G=SO(1,1)\times\R^2$. In the dual $\g^\star$ of
Lie algebra $\g$ of $G$, the orbit $M$ sits as hyperbolic cylinder. In
this picture, the geodesics of the canonical symplectic symmetric
connection $\nabla$ are planar sections of $M\subset\R^3$. The affine
manifold $(M,\nabla)$ turns out to be strictly geodesically convex.

\noindent Moreover, given three points $x,y,z$ in $M$, the equation
$s_xs_ys_z(t)=t$ has always a (unique) solution $t\in M$. In particular,
Weinstein's function $S_{\mathrm{can}}$ is everywhere defined on
$M\times M\times M$. Within suitable global Darboux coordinates
$(M,\omega)\simeq(\R^2,da\wedge d\ell)$, the symmetry at $(a,\ell)$ has
the following expression:
\begin{equation}
s_{(a,\ell)}(a',\ell')=(2a-a',2\cosh(a-a')\ell-\ell')\;;
\end{equation}
while Weinstein's function is given by
\begin{equation}
S_{\mathrm{can}}((a_1,\ell_1),(a_2,\ell_2),(a_3,\ell_3))=
\cyclic_{1,2,3}\sinh(a_1-a_2)\ell_3\;.
\end{equation}
Observe that in this coordinate system one sees how far $M$ is from
being flat: indeed replacing, in the expression of the symmetry map, the
function
\begin{equation}
A^0:M\times M\to\R:(\;(a,\ell)\;,\;(a',\ell')\;)\mapsto\cosh(a-a')
\end{equation}
by the constant function $1$ would yield the flat plane. This function
turns out to be exactly the one which twists the volume form
$\omega\wedge\omega$ on $M\times M$ in the expression of a WKB
quantisation product . More precisely, one has

\begin{theorem}\cite{Biel2} 
There exist Fr\'echet function spaces $\{{\cal A}_\theta\}_{\theta\in\R}$
such that

\begin{enumerate}
\item[(i)] for all $\theta$, one has 
$$
C^\infty_c(M)\subset{\cal A}_\theta\subset C_\infty(M)
$$
\item[(ii)] the formula
$$
(u\star_\theta v)(x)=\frac{1}{\theta^{2}}\int_{M\times M}A^0(y,z)\, \exp\left( 
\frac{i}{\theta}\,S_{\mathrm{can}}(x,y,z)\right)\, u(y)\, v(z) \, dy\, dz
$$
defines an associative product on ${\cal A}_\theta\quad\theta\neq 0$. 
Each pair $({\cal A}_\theta,\star_\theta)$ is then a Fr\'echet algebra.
\item[(iii)] For $u$ and $v$ smooth compactly supported functions on
$M$, one has an asymptotic expansion in powers of $\theta$:
$$
u\star_\theta v\sim 
uv+\frac{\theta}{2i}\{u,v\}\;+\quad\mbox{ \rm higher order terms}
$$
where $\{\;,\;\}$ denotes the Poisson structure associated to the
symplectic form $\omega$.
\end{enumerate}
\end{theorem}
An analysis of the $G$-equivariant (formal) equivalences yields a WKB
expression for every star product on $M$ from the data of the preceding
one. More precisely, one has
\begin{prop}\cite{BDRS2}
Let ${\cal P}\in C^\infty(\R)[[\theta]]$ be a formal function on $\R$.
Then, an asymptotic expansion in powers of $\theta$ of the following
oscillatory integral
{\small
\begin{equation}
(u\star^{\cal P}_\theta v)(x)=\frac{1}{\theta^{2}}\int_{M\times
M}\frac{{\cal P}(a_{\bf x} - a_{\bf z}) {\cal P}(a_{\bf y}-a_{\bf
x})}{{\cal P}(a_{\bf y} -a_{\bf z})} \,A^0(y,z)\, \exp\left(
\frac{i}{\theta}\,S_{\mathrm{can}}(x,y,z)\right)\, u(y)\, v(z) \, dy\,
dz
\end{equation}
}
yields  a $G$-invariant star product on $M$. Moreover, every
$G$-invariant star product on $M$ may be described this way. The choice
\begin{equation}
{\cal P}(a)=\sqrt{\cosh(a)}
\end{equation}
yields a strongly closed star product $\star^{\mbox{\tiny{s.c}}}$ on $M$
(i.e. $\int u\star^{\mbox{\tiny{s.c}}}v=\int u v$).
\end{prop}
Now intrinsically, the amplitude of the oscillating kernel defining the
above strongly closed star product may be described geometrically as
follows. Let
\begin{equation}
\Phi:M\times M\times M\to M\times M\times M:(x,y,z)\mapsto(X,Y,Z)
\end{equation}
with 
\begin{equation}
s_xs_ys_z(X)=X\;,\quad Y=s_z(X)\;,\quad Z=s_y(Y).
\end{equation}

\noindent Now consider the Jacobian map of $\Phi$:
\begin{equation}
{\mbox{\rm
Jac}}_\Phi(x,y,z):=\left|\frac{\partial(X,Y,Z)}{\partial(x,y,z)}\right|.
\end{equation}
Then, it turns out that in the above coordinate system the function
$\mbox{Jac}_\Phi$ depends only on the $a$-coordinates of the points and
that its square root coincides with the above mentioned amplitude. More
precisely, one has
\begin{theorem}\cite{BDRS2}
Let $\theta>0$.
For $u$ and $v$ compactly supported functions on $M$, the formula
\begin{equation}
u\star^{\mbox{\tiny{s.c}}}_\theta v:=\frac{1}{\theta^2}\int_{M\times M}\;\sqrt{{\mbox{\rm
Jac}}_\Phi}\;e^{\frac{i}{\theta}S_{\mathrm{can}}} \;u\otimes v
\end{equation}
extends to $L^2(M)$ as an associative product. The function algebra
$(L^2(M), \star^{\mbox{\tiny{s.c}}}_\theta)$ is then a Hilbert algebra
with respect to the natural Hilbert space structure on $L^2(M)$.
\end{theorem}
From this a continuous field of $C^\star$-algebras deforming $C_0(M)$
may be obtained via standard techniques (see e.g. \cite{Rieffel} for the
flat case and \cite{Biel2} for curved solvable symmetric spaces).

\end{document}